\documentclass{kapedbk} 
\usepackage{amsfonts}
\usepackage{amsmath, amsthm}
\usepackage{algorithm}
\usepackage{algorithmic}

\usepackage[dvips]{graphicx}

\upperandlowercase

\normallatexbib 

\newtheorem{Problem}{Problem}
\newtheorem{Definition}{Definition}
\newtheorem{Theorem}{Theorem}[section]
\newtheorem{Lemma}{Lemma}[section]

\newcommand{\HDiv}[2]{\frac{\left[#1\right]}{\overline{\left[#2\right]}}}

\newcommand{\innprod}[2]{\left <#1, #2 \right>}

\DeclareMathOperator*{\argmax}{argmax}
\DeclareMathOperator*{\argmin}{argmin}

\begin{document}

\articletitle[Descent methods for Nonnegative Matrix Factorization]{Descent methods for Nonnegative Matrix Factorization}

\author{Ngoc-Diep Ho, Paul Van Dooren and Vincent D. Blondel} 
\affil{CESAME, Universit{\'e} catholique de Louvain, \\ Av. Georges Lema\^itre 4, 
B-1348 Louvain-la-Neuve, Belgium. \\ Tel : +32(10) 47 22 64 \ \ Fax : +32(10) 47 21 80}
\email{\{ngoc.ho, paul.vandooren, vincent.blondel\}@uclouvain.be}

\begin{abstract}
In this paper, we present several descent methods that can be applied to nonnegative matrix factorization and we analyze a recently developped fast block coordinate method called Rank-one Residue Iteration (RRI). We also give a comparison of these different methods and show that the new block coordinate method has better properties in terms of approximation error and complexity. By interpreting this method as a rank-one approximation of the residue matrix, we prove that it \emph{converges} and also extend it to the nonnegative tensor factorization and introduce some variants of the method by imposing some additional controllable constraints such as: sparsity, discreteness and smoothness.
\end{abstract}


\section{Introduction}
Linear algebra has become a key tool in almost all modern techniques for data analysis. Most of these techniques make use of linear subspaces represented by eigenvectors of a particular matrix. In this paper, we consider a set of $n$ data points $a_1, a_2, \ldots, a_n$, where each point is a real vector of size $m$, $a_i \in \mathbb{R}^m$. We then approximate these data points by linear combinations of $r$ basis vectors $u_i \in \mathbb{R}^m$:
$$
a_i \approx \sum_{j=1}^r{v_{ij}u_j}, \qquad v_{ij} \in \mathbb{R}, \ u_j\in\mathbb{R}^m.
$$
This can be rewritten in matrix form as $A \approx UV^T$, where $a_i$ and $u_i$ are respectively the columns of $A$ and $U$ and the ${v_{ij}}'s$ are the elements of $V$. Optimal solutions of this approximation in terms of the Euclidean (or Frobenius) norm can be obtained by the Singular Value Decomposition (SVD) \cite{golvanmat}.

In many cases, data points are constrained to a subset of $\mathbb{R}^m$. For example, light intensities, concentrations of substances, absolute temperatures are, by their nature, nonnegative (or even positive) and lie in the nonnegative orthant $\mathbb{R}^m_+$. The input matrix $A$ then becomes elementwise nonnegative and it is then natural to constrain the basis vectors $v_i$ and the coefficients $v_{ij}$ to be nonnegative as well. In order to satisfy this constraint, we need to approximate the columns of $A$ by the following additive model:
$$
a_i \approx \sum_{j=1}^r{v_{ij}u_j}, \qquad v_{ij} \in \mathbb{R}_+, \ u_j\in\mathbb{R}_+^m.
$$
where the $v_{ij}$ coefficients and $u_j$ vectors are nonnegative, $v_{ij} \in \mathbb{R}_+, \ u_j\in\mathbb{R}^m_+$.

Many algorithms have been proposed to find such a representation, which is referred to as a Nonnegative Matrix Factorization (NMF). The earliest algorithms were introduced by Paatero \cite{paatero1994pmf, Paatero1997}. But the topic became quite popular with the publication of the algorithm of Lee and Seung in 1999 \cite{leeseung99nature} where multiplicative rules were introduced to solve the problem. This algorithm is very simple and elegant but it lacks a complete convergence analysis. Other methods and variants can be found in \cite{cjl05b}, \cite{cjl05c}, \cite{hoyer2004}.

The quality of the approximation is often measured by a distance. Two popular choices are the Euclidean (Frobenius) norm and the generalized Kullback-Leibler divergence. In this paper, we focus on the Euclidean distance and we  investigate descent methods for this measure. One characteristic of descent methods is their monotonic decrease until they reach a stationary point. This point maybe located in the interior of the nonnegative orthant or on its boundary. In the second case, the constraints become active and may prohibit any further decrease of the distance measure. This is a key issue to be analyzed for any descent method.

In this paper, $\mathbb{R}_+^m$ denotes the set of nonnegative real vectors (elementwise) and $[v]_+$  the projection of the vector $v$ on $\mathbb{R}_+^m$. We use $v \ge 0$ and $A \ge 0$ to denote nonnegative vectors and matrices and $v > 0$ and $A > 0$ to denote positive vectors and matrices. $A \circ B$ and $\HDiv{A}{B}$ are respectively the Hadamard (elementwise) product and quotient. $A_{:i}$ and $A_{i:}$ are the $i^{th}$ column and $i^{th}$ row of $A$.

This paper is an extension of the internal report \cite{ho2007rep}, where we proposed to decouple the problem based on rank one approximations to create a new algorithm called Rank-one Residue Iteration (RRI). During the revision of this report, we were informed that essentially the same algorithm was independently proposed and published in \cite{Cichocki2007HALS} under the name Hierarchical Alternative Least Squares (HALS). But the present paper gives several additional results wherein the major contributions are \emph{the convergence proof} of the method and its \emph{extensions} to many pratical situations and constraints. The paper also compares a selection of some recent descent methods from the literature and aims at providing a survey of such methods for nonnegative matrix factorizations. For that reason, we try to be self-contained and hence recall some well-known results. We also provide short proofs when useful for a better understanding of the rest of the paper.

We first give a short introduction of low rank approximations, both unconstrained and constrained. In Section \ref{Se:3} we discuss error bounds of various approximations and in Section \ref{Se:4} we give a number of descent methods for Nonnegative Matrix Factorizations. In Section \ref{Se:5} we describe the method based on successive rank one approximations. This method is then also extended to approximate higher order tensor and to take into account other constraints than nonnegativity. In Section \ref{regularize} we discuss various regularization methods and in Section \ref{experiment}, we present numerical experiments comparing the different methods. We end with some concluding remarks.

\section{Low-rank matrix approximation} \label{Lowrank}

Low-rank approximation is a special case of matrix nearness problem \cite{higham1989mnp}. When only a rank constraint is imposed, the optimal approximation with respect to the Frobenius norm can be obtained from the Singular Value Decomposition. 

We first investigate the problem without the nonnegativity constraint on the low-rank approximation. This is useful for understanding properties
of the approximation when the nonnegativity constraints are imposed but inactive. We begin with the well-known Eckart-Young Theorem.

\begin{Theorem}[Eckart-Young] \label{th:EckartYoung} Let $A \in \mathbb{R}^{m \times n}$ ($m \ge n$) have the singular value decomposition
$$
A = P\Sigma Q^T, \ \Sigma = \left( \begin{array}{cccc}
\sigma_1 & 0 & \ldots & 0 \\
0 & \sigma_2 & \ldots & 0 \\
\vdots & \vdots & \ddots & \vdots \\
0 & 0 & \ldots & \sigma_n \\
\vdots & \vdots &  & \vdots \\
0 & 0 & \ldots & 0 
\end{array} \right ) 
$$
where $\sigma_1 \ge \sigma_2 \ge \ldots \ge \sigma_n \ge 0$ are the singular values of $A$ and where $P \in \mathbb{R}^{m \times m}$ and $Q \in \mathbb{R}^{n \times n}$ are orthogonal matrices. Then for $1 \le r \le n$, the matrix
$$
A_r = P\Sigma_r Q^T, \ \Sigma_r = \left( \begin{array}{cccccc}
\sigma_1 & 0 & \ldots & 0 & \ldots & 0\\
0 & \sigma_2 & \ldots & 0 & \ldots & 0 \\
\vdots & \vdots & \ddots & \vdots &  & \vdots\\
0 & 0 & \ldots & \sigma_r & \ldots & 0 \\
\vdots & \vdots &  & \vdots & \ddots & \vdots \\
0 & 0 & \ldots & 0 & \ldots & 0 
\end{array} \right )
$$
is a global minimizer of the problem
\begin{equation} \label{LRApproxProb1}
\min_{B \in \mathbb{R}^{m \times n} \ rank(B) \le r} \frac{1}{2}\|A - B\|_F^2 
\end{equation}
and its error is
$$
\frac{1}{2}\|A - B\|_F^2 = \frac{1}{2}\sum_{i=r+1}^n{\sigma_i^2}.
$$
Moreover, if $\sigma_r > \sigma_{r+1}$ then $A_r$ is the unique global minimizer.
\end{Theorem}
The proof and other implications can be found for instance in \cite{golvanmat}. The columns of $P$ and $Q$ are called singular vectors of $A$, in which vectors corresponding to the largest singular values are referred to as the dominant singular vectors. 

Let us now look at the following modified problem
\begin{equation} \label{LRApproxProb2}
\min_{X \in \mathbb{R}^{m \times r} \ Y\in \mathbb{R}^{n \times r}} \frac{1}{2}\|A - XY^T\|_F^2, 
\end{equation}
where the rank constraint is implicit in the product $XY^T$ since the dimensions of $X$ and $Y$ guarantee that $rank(XY^T) \le r$. Conversely, every matrix of rank less than $r$ can be trivially rewritten as a product $XY^T$, where $X \in \mathbb{R}^{m \times r}$ and $Y\in \mathbb{R}^{n \times r}$. Therefore Problems (\ref{LRApproxProb1}) and (\ref{LRApproxProb2}) are equivalent. But even when the product $A_r = XY^T$ is unique, the pairs $(XR^T, YR^{-1})$ with $R$ invertible, yield the same product $XY^T$. In order to avoid this, we can always choose $X$ and $Y$ such that
\begin{equation} \label{SVDCritical}
X = PD^{\frac{1}{2}} \ \text{and} \ Y = QD^{\frac{1}{2}}, 
\end{equation}
where $P^TP = I_{r \times r}$, $Q^TQ = I_{r \times r}$ and $D$ is $r \times r$ nonnegative diagonal matrix. Doing this is equivalent to computing a compact SVD decomposition of the product $A_r = XY^T = PDQ^T$. 

As usual for optimization problems, we calculate the gradient with respect to $X$ and $Y$ and set them equal to $0$.
\begin{equation} \label{nullgradient}
\nabla_X = XY^TY - AY = 0 \qquad \nabla_Y = YX^TX - A^TX = 0.  
\end{equation}
If we then premultiply $A^T$ with $\nabla_X$  and $A$ with $\nabla_Y$, we obtain
\begin{equation} \label{invarspaces1}
(A^TA)Y = (A^TX)Y^TY \qquad (AA^T)X = (AY)X^TX. 
\end{equation}
Replacing $A^TX = YX^TX$ and $AY = XY^TY$ into (\ref{invarspaces1}) yields
\begin{equation} \label{invarspaces}
(A^TA)Y = YX^TXY^TY \qquad (AA^T)X = XY^TYX^TX. 
\end{equation}
Replacing (\ref{SVDCritical}) into (\ref{invarspaces}) yields
$$
(A^TA)QD^{\frac{1}{2}} = QDP^TPDQ^TQD^{\frac{1}{2}} 
\ \text{and} \ (AA^T)PD^{\frac{1}{2}} = PDQ^TQDP^TPD^{\frac{1}{2}}.
$$
When $D$ is invertible, this finally yields
$$
	(A^TA)Q = QD^2 \ \text{and} \ (AA^T)P = PD^2.
$$

This shows that the columns of $P$ and $Q$ are singular vectors and ${D_{ii}}'s$ are nonzero singular values of $A$.
Notice that if $D$ is singular, one can throw away the corresponding columns of $P$ and $Q$ and reduce it to a smaller-rank approximation with the same properties. Without loss of generality, we therefore can focus on approximations  of Problem (\ref{LRApproxProb2}) which are of exact rank $r$. We can summarize the above reasoning in the following theorem.
\begin{Theorem} \label{th:staLRApprox}
Let $A \in \mathbb{R}^{m \times n}$ ($m > n$ and $rank(A) = t$). If $A_r$ ($1 \le r \le t $) is a rank $r$ stationary point of Problem \ref{LRApproxProb2}, then there exists two orthogonal matrices $P \in \mathbb{R}^{m \times m}$ and $Q \in \mathbb{R}^{n \times n}$ such that:
$$
A = P\hat\Sigma Q^T \ \text{and} \ A_r = P\hat\Sigma_rQ^T
$$
where 
$$\hat\Sigma = \left( \begin{array}{cccc}
\hat\sigma_1 & 0 & \ldots & 0 \\
0 & \hat\sigma_2 & \ldots & 0 \\
\vdots & \vdots & \ddots & \vdots \\
0 & 0 & \ldots & \hat\sigma_n \\
\vdots & \vdots &  & \vdots \\
0 & 0 & \ldots & 0 
\end{array} \right ) 
, \qquad \hat\Sigma_r = \left( \begin{array}{cccccc}
\hat\sigma_1 & 0 & \ldots & 0 & \ldots & 0\\
0 & \hat\sigma_2 & \ldots & 0 & \ldots & 0 \\
\vdots & \vdots & \ddots & \vdots &  & \vdots\\
0 & 0 & \ldots & \hat\sigma_r & \ldots & 0 \\
\vdots & \vdots &  & \vdots & \ddots & \vdots \\
0 & 0 & \ldots & 0 & \ldots & 0 
\end{array} \right )
$$
and the ${\hat\sigma_i}'s$ are unsorted singular values of $A$. Moreover, the approximation error is:
$$
\frac{1}{2}\|A - A_r\|_F^2 = \frac{1}{2}\sum_{i=r+1}^t{\hat\sigma_i^2}.
$$
\end{Theorem}
This result shows that, if the singular values are all different, there are $\frac{n!}{r!(n-r)!}$ possible stationary points $A_r$. When there are multiple singular values, there will be infinitely many stationary points $A_r$ since there are infinitely many singular subspaces. The next result will identify the minima among all stationary points. Other stationary points are saddle points whose every neighborhood contains both smaller and higher points.
\begin{Theorem} \label{globalmin}
The only minima of Problem \ref{LRApproxProb2} are given by Theorem \ref{th:EckartYoung} and are global minima. All other stationary points are saddle points.
\end{Theorem}

\begin{proof} 
Let us assume that $A_r$ is a stationary point given by Theorem \ref{th:staLRApprox} but not by Theorem \ref{th:EckartYoung}. Then there always exists a permutation of the columns of $P$ and $Q$, and of the diagonal elements
of $\hat\Sigma$ and $\hat\Sigma_r$ such that $\hat\sigma_{r+1} > \hat\sigma_r$. We then construct two points in the $\epsilon$-neighborhood of $A_r$ that yield an increase and a decrease, respectively, of the distance measure. They are obtained by taking:
$$
\overline{\Sigma}_r(\epsilon) = \left( \begin{array}{ccccc}
\hat\sigma_1 + \epsilon  & \ldots  & 0 & \ldots & 0\\
\vdots  & \ddots & \vdots & \ldots & \vdots \\
0 & \ldots &  \hat\sigma_r & \ldots & 0 \\
\vdots &  \vdots & \vdots & \ddots & \vdots\\
0 & \ldots & 0 & \ldots & 0 
\end{array} \right), \qquad \overline{A}_r(\epsilon) = P\overline{\Sigma}_r(\epsilon)Q^T
$$
and 
$$
\underline{\Sigma}_r(\epsilon) = \left( \begin{array}{cccccc}
\hat\sigma_1 & \ldots & 0 & 0 & \ldots & 0\\
\vdots & \ddots & \vdots & \vdots & \ldots & \vdots \\
0  & \ldots & \hat\sigma_r & \epsilon \sqrt{\hat\sigma_r}& \vdots & 0\\
0 &  \ldots & \epsilon \sqrt{\hat\sigma_r} & \epsilon^2 & \ldots & 0 \\
\vdots & \vdots & \vdots & \vdots & \ddots & \vdots\\
0 & 0 & \ldots & 0 & \ldots & 0 
\end{array} \right ), \qquad \underline{A}_r(\epsilon) = P\underline{\Sigma}_r(\epsilon)Q^T.
$$
Clearly $\overline{A}_r(\epsilon)$ and $\underline{A}_r(\epsilon)$ are of rank $r$. Evaluating the distance measure yields
\begin{eqnarray*}
\|A - \underline{A}_r(\epsilon)\|_F^2 &=& 2 \hat\sigma_r \epsilon^2 + (\hat\sigma_{r+1} - \epsilon^2)^2 + \sum_{i = r+2}^t\hat\sigma_i^2 \\ 
& = & \epsilon^2[\epsilon^2-2(\hat\sigma_{r+1} - \hat\sigma_{r})] + \sum_{i = r+1}^t\hat\sigma_i^2 \\
& < & \sum_{i = r+1}^t\hat\sigma_i^2  = \|A - A_r\|_F^2\\
\end{eqnarray*}
for all $\epsilon \in (0, \sqrt{2(\hat\sigma_{r+1} - \hat\sigma_{r})})$ and 
$$
\|A - \overline{A}_r(\epsilon)\|_F^2 = \epsilon^2 + \sum_{i = r+1}^t\hat\sigma_i^2 > \sum_{i = r+1}^t\hat\sigma_i^2 = \|A - A_r\|_F^2
$$
for all $\epsilon > 0$. Hence, for an arbitrarily small positive $\epsilon$, we obtain 
$$
\|A - \underline{A}_r(\epsilon)\|_F^2 < \|A - A_r\|_F^2 < \|A - \overline{A}_r(\epsilon)\|_F^2
$$
which shows that $A_r$ is a saddle point of the distance measure. 
\end{proof}

When we add a nonnegativity constraint in the next section, the results of this section will help to identify stationary points at which all the nonnegativity constraints are inactive.

\section{Nonnegativity constraint} \label{Se:3}

In this section, we investigate the problem of Nonnegative Matrix Factorization. This problem differs Problem \ref{LRApproxProb2} in the previous section because of the additional nonnegativity constraints on the factors. We
first discuss the effects of adding such a constraint. By doing so, the problem is no longer easy because of the existence of local minima at the boundary of the nonnegative orthant. Determining the lowest minimum among these minima is far from trivial. On the other hand, a minimum that coincides with a minimum of the unconstrained problem (i.e. Problem \ref{LRApproxProb2}) may be easily reached by standard descent methods, as we will see.

\begin{Problem}[Nonnegative matrix factorization - NMF] Given a $m \times n$ nonnegative matrix $A$ and an integer $r < \min(m, n)$, solve
$$
\min_{U \in \mathbb{R}^{m \times r}_+ \ V \in \mathbb{R}^{n \times r}_+}{\frac{1}{2}\|A - UV^T\|^2_F} \label{prob:NMF}.
$$
\end{Problem}

Where $r$ is called the reduced rank. From now on, $m$ and $n$ will be used to denote the size of the target matrix $A$ and $r$ is the reduced rank of a factorization. 

We rewrite the nonnegative matrix factorization as a standard nonlinear optimization problem:
$$
\min_{-U \le 0 \ -V \le 0}{\frac{1}{2}\|A - UV^T\|^2_F}.
$$
The associated Lagrangian function is
$$
L(U, V, \mu, \nu) = {\frac{1}{2}\|A - UV^T\|^2_F} - \mu \circ U - \nu \circ V,
$$
where $\mu$ and $\nu$ are two \emph{matrices} of the same size of $U$ and $V$, respectively, containing the Lagrange multipliers associated with the nonnegativity constraints $U_{ij} \ge 0$ and $V_{ij} \ge 0$. Then the Karush-Kuhn-Tucker conditions for the nonnegative matrix factorization problem say that if $(U, V)$ is a local minimum, then there exist $\mu_{ij} \ge 0$ and $\nu_{ij} \ge 0$ such that:
\begin{eqnarray}
U \ge 0&, & \quad V \ge 0,  \label{ch2-1}\\
\nabla L_U = 0&, & \quad \nabla L_V = 0,  \label{ch2-2}\\
\mu \circ U = 0&,  & \quad \nu \circ V  = 0 \label{ch2-3}.
\end{eqnarray}
Developing (\ref{ch2-2}) we have:
$$
AV - UV^TV - \mu = 0,  \quad A^TU - VU^TU - \nu = 0\\
$$
or 
$$
\mu = -(UV^TV - AV),  \quad \nu = -(VU^TU - A^TU). \\
$$
Combining this with $\mu_{ij} \ge 0$, $\nu_{ij} \ge 0$ and (\ref{ch2-3}) gives the following conditions:
\begin{eqnarray}
U \ge 0&, & \quad V \ge 0, \label{KKT1} \\
\nabla F_U = UV^TV - AV \ge 0&, & \quad \nabla F_V = VU^TU - A^TU \ge 0, \label{KKT2} \\
U\circ(UV^TV - AV) = 0&,  & \quad V\circ (VU^TU - A^TU)  = 0, \label{KKT3}
\end{eqnarray}
where the corresponding Lagrange multipliers for $U$ and $V$ are also the gradient of $F$ with respect to $U$ and $V$.
Since the Euclidean distance is not convex with respect to both variables $U$ and $V$ at the same time, these conditions are only necessary. This is implied because of the existence of saddle points and maxima. We then call all the points that satisfy the above conditions, the stationary points.

\begin{Definition}[NMF stationary point]
We call $(U, V)$ a stationary point of the NMF Problem if and only if $U$ and $V$ satisfy the KKT conditions (\ref{KKT1}), (\ref{KKT2}) and (\ref{KKT3}).
\end{Definition}

Alternatively, a stationary point $(U, V)$ of the NMF problem can also be defined by using the following necessary condition (see for example \cite{bertsekas1999np}) on the convex sets $\mathbb{R}_+^{m \times r}$ and $\mathbb{R}_+^{n \times r}$, that is
\begin{equation}
\innprod{\left (\begin{array}{c} \nabla F_{U} \\ \nabla F_{V}\end{array} \right)}{\left (\begin{array}{c} X - U \\ Y - V\end{array} \right)} \ge 0, \qquad \forall \ X \in \mathbb{R}_+^{m \times r}, \ Y \in \mathbb{R}_+^{n \times r}, \label{fonmfcond}
\end{equation}
which can be shown to be equivalent to the KKT conditions (\ref{KKT1}), (\ref{KKT2}) and (\ref{KKT3}). Indeed, it is trivial that the KKT conditions imply (\ref{fonmfcond}). And by carefully choosing different values of $X$ and $Y$ from (\ref{fonmfcond}), one can easily prove that the KKT conditions hold.

There are two values of reduced rank $r$ for which we can trivially identify the global solution which are $r = 1$ and $r = \min(m, n)$. For $r = 1$, a pair of dominant singular vectors are a global minimizer. And for $r = \min(m, n)$, $(U = A, V = I)$ is a global minimizer. Since most of existing methods for the nonnegative matrix factorization are descent algorithms, we should pay attention to all local minimizers. For the rank-one case, they can easily be characterized.

\subsection{Rank one case}

The rank-one NMF problem of a nonnegative matrix $A$ can be rewritten as
\begin{equation}
\min_{u \in \mathbb{R}_+^m \ v \in \mathbb{R}_+^n}{\frac{1}{2}}\|A - uv^T\|_F^2 \label{R1NMF}
\end{equation}
and a complete analysis can be carried out. It is well known that any pair of nonnegative Perron vectors of $AA^T$ and $A^TA$ yields a global minimizer of this problem, but we can also show that the \emph{only} stationary points of (\ref{R1NMF}) are given by such vectors. The following theorem excludes the case where $u=0$ and/or $v=0$.

\begin{Theorem} \label{th:R1NMF}
The pair $(u, v)$ is a local minimizer of (\ref{R1NMF}) if and only if $u$ and $v$ are nonnegative eigenvectors of $AA^T$ and $A^TA$ respectively of the eigenvalue $\sigma=\|u\|_2^2\|v\|_2^2$.
\end{Theorem}

\begin{proof}
The \emph{if part} easily follows from Theorem \ref{th:staLRApprox}. For the \emph{only if part} we proceed as follows. Without loss of generality, we can permute the rows and columns of $A$ such that the corresponding vectors $u$ and $v$ are partitioned as $(u_+ \ 0)^T$ and $(v_+ \ 0)^T$ respectively, where $u_+$, $v_+ > 0$. Partition the corresponding matrix $A$ conformably as follows
$$
A = \left ( \begin{array}{cc}A_{11} & A_{12} \\ A_{21} & A_{22} \end{array} \right ),
$$
then from (\ref{KKT2}) we have
$$
\left ( \begin{array}{cc} u_+v_+^T & 0 \\ 0 & 0 \end{array}\right ) \left ( \begin{array}{c} v_+ \\ 0 \end{array}\right ) - \left ( \begin{array}{cc}A_{11} & A_{12} \\ A_{21} & A_{22} \end{array} \right ) \left ( \begin{array}{c} v_+ \\ 0 \end{array}\right ) \ge 0
$$
and 
$$
\left ( \begin{array}{cc} v_+u_+^T & 0 \\ 0 & 0 \end{array}\right ) \left ( \begin{array}{c} u_+ \\ 0 \end{array}\right ) - \left ( \begin{array}{cc}A^T_{11} & A^T_{21} \\ A_{12}^T & A^T_{22} \end{array} \right ) \left ( \begin{array}{c} u_+ \\ 0 \end{array}\right ) \ge 0
$$

implying that  $A_{21}v_+ \le 0$ and $A^T_{12}u_+ \le 0$. Since $A_{21}$ , $A_{12} \ge 0$ and $u_+$, $v_+ > 0$, we can conclude that $A_{12} = 0$ and $A_{21} = 0$. Then from (\ref{KKT3}) we have:
$$
u_+ \circ (\|v_+\|_2^2u_+ - A_{11}v_+) = 0 \ \text{and} \ v_+ \circ (\|u_+\|_2^2v_+ - A^+_{11}u_+) = 0.
$$
Since $u_+$, $v_+ > 0$, we have:
$$
\|v_+\|_2^2u_+ = A_{11}v_+ \ \text{and} \ \|u_+\|_2^2v_+ = A_{11}^Tu_+
$$
or
$$
\|u_+\|_2^2\|v_+\|_2^2u_+ = A_{11}A_{11}^Tu_+ \ \text{and} \ \|u_+\|_2^2\|v_+\|_2^2v_+ = A_{11}^TA_{11}v_+.
$$
Setting $\sigma = \|u_+\|_2^2\|v_+\|_2^2$ and using the block \emph{diagonal} structure of $A$ yields the desired result. 
\end{proof}

Theorem \ref{th:R1NMF} guarantees that all stationary points of the rank-one case are nonnegative singular vectors of a submatrix of $A$. These results imply that a global minimizer of the rank-one NMF can be calculated correctly based on the largest singular value and corresponding singular vectors of the matrix $A$.

For ranks other than $1$ and $\min(m, n)$, there are no longer trivial stationary points. In the next section, we try to derive some simple characteristics of the local minima of the nonnegative matrix factorization.

The KKT conditions (\ref{KKT3}) help to characterize the stationary points of the NMF problem. Summing up all the elements of one of the conditions (\ref{KKT3}), we get:
\begin{eqnarray}
\nonumber 0 &=& \sum_{ij} \left (U\circ(UV^TV - AV) \right)_{ij} \\
\nonumber   &=& \left<U, UV^TV - AV \right>\\
                    &=& \left<UV^T, UV^T - A \right> \label{sumkkt}.
\end{eqnarray}

From that, we have some simple characteristics of the NMF solutions:

\begin{Theorem} \label{th:solball}
Let $(U, V)$ be a stationary point of the NMF problem, then $UV^T \in \mathcal{B}\left (\frac{A}{2}, \frac{1}{2}\|A\|_F \right)$, the ball centered at $\frac{A}{2}$ and with radius = $\frac{1}{2}\|A\|_F$.
\end{Theorem}
\proof{
From (\ref{sumkkt}) it immediately follows that
$$
\left<\frac{A}{2} - UV^T, \frac{A}{2} - UV^T\right> = \left< \frac{A}{2}, \frac{A}{2}\right>
$$
which implies 
$$
UV^T \in \mathcal{B}\left (\frac{A}{2}, \frac{1}{2}\|A\|_F \right).
$$ 
} 

\begin{Theorem} \label{th:solnorm}
Let $(U, V)$ be a stationary of the NMF problem, then
$$
\frac{1}{2} \|A - UV^T\|_F^2 = \frac{1}{2}(\|A\|_F^2 - \|UV^T\|_F^2).
$$
\end{Theorem}
\begin{proof}
From (\ref{sumkkt}), we have $\left<UV^T, A \right> = \left<UV^T, UV^T \right>$. Therefore,
\begin{eqnarray*}
\frac{1}{2}\left<A - UV^T, A - UV^T\right> &=& \frac{1}{2}(\|A\|_F^2 - 2\left<UV^T, A\right> + \|UV^T\|_F^2) \\ &=& \frac{1}{2}(\|A\|_F^2 - \|UV^T\|_F^2).
\end{eqnarray*} 
\end{proof}

Theorem \ref{th:solnorm} also suggests that at a stationary point $(U, V)$ of the NMF problem, we should have $\|A\|_F^2 \ge \|UV^T\|_F^2$. This norm inequality can be also found in \cite{catral2004rrn} for less general cases where we have $\nabla F_U = 0$ and $\nabla F_V = 0$ at a stationary point. For this particular class of NMF stationary point, all the nonnegativity constraints on $U$ and $V$are inactive. And all such stationary points are also stationary points of the unconstrained problem, characterized by Theorem \ref{th:staLRApprox}.

We have seen in Theorem \ref{th:staLRApprox} that, for the unconstrained least-square problem the only stable stationary points are in fact global minima. Therefore, if the stationary points of the constrained problem are inside the nonnegative orthant (i.e. all constraints are inactive), we can then probably reach the global minimum of the NMF problem. This can be expected because the constraints may no longer prohibit the descent of the update.

Let $A_r$ be the optimal rank-$r$ approximation of a nonnegative matrix $A$, which we obtain from the singular value decomposition, as indicated in Theorem \ref{th:staLRApprox}. Then we can easily construct its nonnegative part $[A_r]_+$, which is obtained from $A_r$ by just setting all its negative elements equal to zero. This is in fact the closest matrix in the cone of nonnegative matrices to the matrix $A_r$, in the Frobenius norm (in that sense, it is its projection on that cone). We now derive some bounds for the error $\|A-[A_r]_+\|_F$.

\begin{Theorem} \label{th:Ar+}
Let $A_r$ be the best rank $r$ approximation of a nonnegative matrix $A$, and let $[A_r]_+$ be its nonnegative part, then
$$
 \|A - [A_r]_+\|_F \leq \|A-A_r\|_F.
$$
\end{Theorem}

\begin{proof}
This follows easily from the convexity of the cone of nonnegative matrices. Since both $A$ and $[A_r]_+$
are nonnegative and since $[A_r]_+$ is the closest matrix in that cone to $A_r$ we immediately obtain
the inequality
$$ \|A - A_r\|_F^2  \geq \|A-[A_r]_+\|_F^2 + \|A_r-[A_r]_+ \|_F^2  \geq \|A-[A_r]_+\|_F^2 $$ 
from which the result readily follows.
\end{proof}

The approximation $[A_r]_+$ has the merit of requiring as much storage as a rank $r$ approximation, even though
its rank is larger than $r$ whenever $A_r\neq [A_r]_+$. We will look at the quality of this approximation in 
Section \ref{experiment}.
If we now compare this bound with the nonnegative approximations then we obtain the following inequalities. Let $U_*V_*^T$ be an optimal nonnegative rank $r$ approximation of $A$ and let $UV^T$ be any stationary point of the KKT conditions for a nonnegative rank $r$ approximation, then we have :
$$  \|A - [A_r]_+\\|_F^2  \leq \|A-A_r\|_F^2 = \sum_{i=r+1}^n\sigma_i^2 \leq \|A-U_*V_*^T \|_F^2  \leq \|A-UV^T\|_F^2.
$$

For more implications of the NMF problem, see \cite{ho2008thesis}. 

\section{Existing descent algorithms} \label{Se:4}

We focus on descent algorithms that guarantee a non increasing update at each iteration. Based on the search space, we have two categories: \textit{Full-space search} and \textit{(Block) Coordinate search}.

Algorithms in the former category try to find updates for both $U$ and $V$ at the same time. This requires a search for a descent direction in the $(m + n)r$-dimensional space. Note also that the NMF problem in this full space is not convex but the optimality conditions may be easier to achieve.

Algorithms in the latter category, on the other hand, find updates for each (block) coordinate in order to guarantee the descent of the objective function. Usually, search subspaces are chosen to make the objective function convex so that efficient methods can be applied. Such a simplification might lead to the loss of some convergence properties. Most of the algorithms use the following column partitioning:
\begin{equation}
\frac{1}{2}\|A - UV^T\|^2_F = \frac{1}{2}\sum_{i=1}^n{\|A_{:,i} - U(V_{i,:})^T\|^2_2}, \label{colpartition}
\end{equation}
which shows that one can minimize with respect to each of the rows of $V$ independently. The problem thus decouples into smaller convex problems. This leads to the solution of quadratic problems of the form
\begin{equation}
\min_{v \ge 0} {\frac{1}{2}\|a - Uv\|^2_2}. \label{smallprob}
\end{equation}
Updates for the rows of $V$ are then alternated with updates for the rows of $U$ in a similar manner by transposing $A$ and $UV^T$.

Independent on the search space, most of algorithms use the Projected Gradient scheme for which three basic steps are carried out in each iteration:
\begin{itemize}
\item Calculating the gradient $\nabla F(x^k)$,
\item Choosing the step size $\alpha^k$,
\item Projecting the update on the nonnegative orthant $$x^{k+1} = [x^k - \alpha^k\nabla F(x^k)]_+,$$
\end{itemize}

where $x^k$ is the variable in the selected search space. The last two steps can be merged in one iterative process and must guarantee a sufficient decrease of the objective function as well as the nonnegativity of the new point.

\subsection{Multiplicative rules (Mult)}

Multiplicative rules were introduced in \cite{leeseung99nature}. The algorithm applies a block coordinate type search and uses the above column partition to formulate the updates. A special feature of this method is that the step size is calculated for each element of the vector. For the elementary problem (\ref{smallprob}) it is given by
$$
v^{k+1} = v^{k} - \mathbf{\alpha}^k \circ \nabla F(v^{k+1}) = v^k \circ \HDiv{U^Ta}{U^TUv^k}
$$
where $[\mathbf{\alpha}^k]_i = \frac{v_i}{[U^TUv]_i}$. Applying this to all rows of $V$ and $U$ gives the updating rule of Algorithm \ref{alg:LeeSeung} to compute
$$(U^*, V^*) = \argmin_{U \ge 0 \ V \ge 0}{\|A - UV^T\|^2_F}.$$

\begin{algorithm}
\caption{(Mult)}
\label{alg:LeeSeung}
\begin{algorithmic}[1]
\STATE Initialize $U^0$, $V^0$ and $k = 0$
\REPEAT
\STATE $U^{k+1} = U^{k} \circ \HDiv{AV^k}{U^{k}(V^k)^T(V^k)}$
\STATE $V^{k+1} = V^{k} \circ \HDiv{A^TU^{k+1}}{V^{k}(U^{k+1})^T(U^{k+1})}$
\STATE $k = k + 1$
\UNTIL{Stopping condition}
\end{algorithmic}
\end{algorithm}

These updates guarantee automatically the nonnegativity of the factors but may fail to give a sufficient decrease of the objective function. It may also get stuck in a non-stationary point and hence suffer from
a poor convergence. Variants can be found in \cite{cjl05c, merritt2005ipg}.

\subsection{Line search using Armijo criterion (Line)}

In order to ensure a sufficient descent, the following projected gradient scheme with Armijo criterion \cite{cjl05b, cjl07} can be applied to minimize
$$x^* = \argmin_{x}{F(x)}.$$

\begin{algorithm}
\caption{(Line)}
\label{alg:ArmijoSearch}
\begin{algorithmic}[1]
\STATE Initialize $x^0$, $\sigma$, $\beta$, $\alpha_0 = 1$ and $k = 1$
\REPEAT
\STATE $\alpha_k = \alpha_{k - 1}$
\STATE $y = [x^k - \alpha_k \nabla F(x^k)]_+$
\IF{$F(y) - F(x^k) > \sigma \left<\nabla F(x^k), y - x^k \right>$}
    \REPEAT
        \STATE $\alpha_k = \alpha_k \cdot \beta$
        \STATE $y = [x^k - \alpha_k \nabla F(x^k)]_+$
    \UNTIL {$F(y) - F(x^k) \le \sigma \left<\nabla F(x^k), y - x^k \right>$}
\ELSE
    \REPEAT
        \STATE $lasty = y$
        \STATE $\alpha_k = \alpha_k / \beta$
        \STATE $y = [x^k - \alpha_k \nabla F(x^k)]_+$
    \UNTIL {$F(y) - F(x^k) > \sigma \left<\nabla F(x^k), y - x^k \right>$}
    \STATE $y = lasty$
\ENDIF
\STATE $x^{k+1} = y$
\STATE $k = k + 1$
\UNTIL{Stopping condition}
\end{algorithmic}
\end{algorithm}

Algorithm \ref{alg:ArmijoSearch} needs two parameters $\sigma$ and
$\beta$ that may affect its convergence. It requires only the
gradient information, and is applied in \cite{cjl05b} for two
different strategies : for the whole space $(U, V)$ (Algorithm FLine) and for $U$ and
$V$ separately in an alternating fashion (Algorithm CLine). With a good choice of
parameters ($\sigma = 0.01$ and $\beta = 0.1$) and a good strategy
of alternating between variables, it was reported in \cite{cjl05b}
to be the faster than the multiplicative rules.

\subsection{Projected gradient with first-order approximation (FO)}

In order to find the solution to
$$x^* = \argmin_{x}{F(x)}$$
we can also approximate at each iteration the function $F(X)$ using:
$$
\tilde F(x) = F(x^k) + \left<\nabla_x F(x^k), x-x^k\right> + \frac{L}{2}\|x^k - x\|_2^2,
$$
where $L$ is a Lipshitz constant satisfying $F(x) \le \tilde{F}(x), \ \forall x$. 
Because of this inequality, the solution of the following problem
$$
x_{k+1} = \argmin_{x \ge 0}{\tilde F(x)}
$$
also is a point of descent for the function $F(x)$ since
$$ F(x_{k+1}) \leq  \tilde F(x_{k+1}) \leq \tilde F(x_{k}) = F(x_{k}).
$$ 
Since the constant $L$ is not known a priori, an inner loop is
needed. Algorithm \ref{FirstOrderApprox} presents an iterative way
to carry out this scheme. As in the previous algorithm this also
requires only the gradient information and can therefore can be
applied to two different strategies: to the whole space $(U, V)$ (Algorithm FFO) and
to $U$ and $V$ separately in an alternating fashion (Algorithm CFO).

\begin{algorithm}
\caption{(FO)}
\label{FirstOrderApprox}
\begin{algorithmic}[1]
\STATE Initialize $x^0$, $L_0$ and $k = 0$
\REPEAT
\STATE $y = [x^k - \frac{1}{L_k} \nabla F(x^k)]_+$
    \WHILE{$F(y) - F(x^k) >  \left<\nabla F(x^k), y - x^k \right> + \frac{L_k}{2}\|y - x^k\|_2^2$}
    \STATE $L_k = L_k / \beta$
    \STATE $Y = [x^k - \frac{1}{L_k} \nabla F(x^k)]_+$
    \ENDWHILE
\STATE $x^{k+1} = y$
\STATE $L_{k+1} = L_k \cdot \beta$
\STATE $k = k + 1$
\UNTIL{Stopping condition}
\end{algorithmic}
\end{algorithm}

A main difference with the previous algorithm is its stopping
criterion for the inner loop.  This algorithm requires also a
parameter $\beta$ for which the practical choice is $2$.

\subsection{Alternative least squares methods}

The first algorithm proposed for solving the nonnegative matrix factorization was the alternative least squares method \cite{paatero1994pmf}. It is known that, fixing either $U$ or $V$, the problem becomes a least squares problem with nonnegativity constraint. 
 
\begin{algorithm}[htb]
\caption{Alternative Least Square (ALS)}
\label{alg:ALS}
\begin{algorithmic}[1]
\STATE Initialize $U$ and $V$
\REPEAT
\STATE Solve: $\min_{V \ge 0}{\frac{1}{2}\|A - UV^T\|_F^2}$ \label{ALS1}
\STATE Solve: $\min_{U \ge 0}{\frac{1}{2}\|A^T - VU^T\|_F^2}$ \label{ALS2}
\UNTIL{Stopping condition}
\end{algorithmic}
\end{algorithm}

Since the least squares problems in Algorithm \ref{alg:ALS} can be perfectly decoupled into smaller problems corresponding to the columns or rows of $A$, we can directly apply methods for the Nonnegative Least Square problem to each of the small problem. Methods that can be applied are \cite{lawson1974sls}, \cite{bro1997fnn}, etc. 

\subsection{Implementation}

The most time-consuming job is the test for the sufficient decrease, which is also the stopping condition for the inner loop. As mentioned at the beginning of the section, the above methods can be carried out using two different strategies: full space search or coordinate search. In some cases, it is required to evaluate repeatedly the function $F(U, V)$. We mention here how to do this efficiently with the coordinate search.

\textbf{Full space search}: The exact evaluation of $F(x) = F(U, V)$ $= \|A - UV^T\|_F^2$ need $O(mnr)$ operations. When there is a correction $y = (U + \Delta U, V + \Delta V)$, we have to calculate $F(y)$ which also requires $O(mnr)$ operations. Hence, it requires $O(tmnr)$ operations to determine a stepsize in $t$ iterations of the inner loop.

\textbf{Coordinate search}: when $V$ is fixed, the Euclidean distance is a quadratic function on $U$: 
\begin{eqnarray*}
F(U) = \|A - UV^T\|_F^2 &=& \innprod{A}{A} - 2 \innprod{UV^T}{A} + \innprod{UV^T}{UV^T} \\
&=& \|A\|_F^2 - 2 \innprod{U}{AV} + \innprod{U}{U(V^TV)}.
\end{eqnarray*}
The most expensive step is the computation of $AV$, which requires $O(mnr)$ operations. But when $V$ is fixed, $AV$ can be calculated once at the beginning of the inner loop. The remaining computations are $\innprod{U}{AV}$ and $\innprod{U}{U(V^TV)}$, which requires $O(nr)$ and $O(nr^2 + nr)$ operations. Therefore, it requires $O(tnr^2)$ operations to determine a stepsize in $t$ iterations of the inner loop which is much less than $O(tmnr)$ operations. This is due to the assumption $r \ll n$. Similarly, when $U$ fixed, $O(tmr^2)$ operations are needed to determine a stepsize.

If we consider an iteration is a sweep, i.e. once all the variables are updated, the following table summarizes the complexity of each sweep of the  described algorithms:

\begin{center}
\begin{tabular}{l|c}
Algorithm & Complexity per iteration\\
\hline
Mult & $O(mnr)$ \\
FLine & $O(tmnr)$ \\
CLine & $O(t_1nr^2 + t_2mr^2)$ \\
FFO & $O(tmnr)$ \\
CFO & $O(t_1nr^2 + t_2mr^2)$ \\
ALS & $O(2^rmnr)^*$ \\
IALS & $O(mnr)$\\
\end{tabular}
\end{center}

where $t$, $t_1$ and $t_2$ are the number of iterations of inner loops, which can not be bounded in general. For algorithm $ALS$, the complexity is reported for the case where the active set method \cite{lawson1974sls} is used. Although $O(2^rmnr)$ is a very high \emph{theorical} upper bound that count all the possible subsets of $r$ variables of each subproblem, in practice, the active set method needs much less iterations to converge. One might as well use more efficient convex optimization tools to solve the subproblems instead of the active set method. 

\subsection{Scaling and Stopping criterion} \label{scalingstopcond}

For descent methods, several stopping conditions are used in the literature.  We now discuss some problems when implementing these conditions for NMF.

The very first condition is the decrease of the objective function. The algorithm should stop when it fails to make the objective function decrease with a certain amount~:
$$
F(U^{k+1}, V^{k+1}) - F(U^k, V^k) < \epsilon \quad or \quad
\frac{F(U^{k+1}, V^{k+1}) - F(U^k, V^k)}{F(U^k, V^k)} < \epsilon.
$$
This is not a good choice for all cases since the algorithm may stop at a point very far from a stationary point. Time and iteration bounds can also be imposed for very slowly converging algorithms. But here again this may not be good for the optimality conditions. A better choice is probably the norm of the projected gradient as suggested in \cite{cjl05b}. For the NMF problem it is defined as follows~:
$$
[\nabla^P_X]_{ij} = \left \{ \begin{array}{ll}
[\nabla_X]_{ij} & \text{if } X_{ij} > 0 \\
\min(0, [\nabla_X]_{ij}) & \text{if } X_{ij} = 0
\end{array} \right.
$$
where $X$ stands for $U$ or $V$. The proposed condition then becomes
\begin{equation}
\left \|\left( \begin{array}{c} \nabla^P_{U^k} \\ \nabla^P_{V^k} \end{array}\right) \right\|_F \le \epsilon \left \|\left( \begin{array}{c} \nabla_{U^1} \\ \nabla_{V^1}\end{array}\right) \right\|_F. 
\label{projnormcond}
\end{equation}
We should also take into account the scaling invariance between $U$
and $V$.  Putting $\bar U = \gamma U$ and $\bar V =
\frac{1}{\gamma}V$ does not change the approximation $UV^T$ but the
above projected gradient norm is affected:
\begin{eqnarray}
\left \|\left (\begin{array}{c}\nabla^P_{\bar U} \\ \nabla^P_{\bar V} \end{array}\right) \right\|^2_F&=&\|\nabla^P_{\bar U}\|^2_F +
\|\nabla^P_{\bar V}\|^2_F = \frac{1}{\gamma^2}\|\nabla^P_{U}\|^2_F +
\gamma^2 \|\nabla^P_{V}\|^2_F \label{newstopcond} \\&\neq&\left \|\left (\begin{array}{c}\nabla^P_{U} \\ \nabla^P_{ V} \end{array}\right) \right\|^2_F.
\nonumber
\end{eqnarray}
Two approximate factorizations $UV^T = \bar U\bar V^T$ resulting in the same approximation should be considered equivalent in terms of precision. One could choose
$\gamma^2:=\|\nabla^P_{U}\|_F/\|\nabla^P_{V}\|_F$, which minimizes
(\ref{newstopcond}) and forces $\|\nabla^P_{\bar U}\|_F=\|\nabla^P_{\bar V}\|_F$, 
but this may not be a good choice when only one of the gradients
$\|\nabla^P_{\bar U}\|_F$ and $\|\nabla^P_{\bar V}\|_F$ is nearly
zero.

In fact, the gradient $\left( \begin{array}{c} \nabla_{U} \\ \nabla_{V}\end{array}\right)$ 
is scale dependent in the NMF problem and any stopping criterion that uses gradient
information is affected by this scaling. To limit that effect, we
suggest the following scaling after each iteration:
$$
\tilde U_k \leftarrow U_k D_k \qquad \tilde V_k \leftarrow V_k
D_k^{-1}
$$
where $D_k$ is a positive diagonal matrix:
$$
[D_k]_{ii} = \sqrt{\frac{\|V_{:i}\|_2}{\|U_{:i}\|_2}}.
$$
This ensures that $\|\tilde U_{:i}\|^2_F=\|\tilde V_{:i}\|^2_F$ and hopefully reduces also the difference between $\|\nabla^P_{\tilde U}\|^2_F$ and $\|\nabla^P_{\tilde V}\|^2_F$. Moreover, it may help to avoid 

The same scaling should be applied to the initial point as well $(U_1, V_1)$ when using (\ref{projnormcond}) as the stopping condition.
 
\section{Rank-one Residue Iteration} \label{Se:5}

In the previous section, we have seen that it is very appealing to decouple the problem into convex subproblems. But this may
``converge'' to solutions that are far from the global minimizers of the problem.

In this section, we analyze a different decoupling of the problem based on rank one approximations. This also allows us to formulate a very simple basic subproblem. This scheme has a major advantage over other methods~: the subproblems can be optimally solved in closed form. Therefore it can be proved to have a strong convergence results through its \emph{damped} version and it can be extended to more general types of factorizations such as for nonnegative tensors and to some practical constraints such as sparsity and smoothness. Moreover, the experiments in Section \ref{experiment} suggest that this method outperforms the other ones in most cases. During the completion of the revised version of this report, we were informed that an independent report \cite{Cichocki2007HALS} had also proposed this decoupling without any convergence investivation and extentions.

\subsection{New partition of variables}

Let the $u_i$'s and $v_i$'s be respectively the columns of $U$ and $V$.
Then the NMF problem can be rewritten as follows~:

\begin{Problem}[Nonnegative Matrix Factorization] Given a $m \times n$ nonnegative matrix $A$, solve
$$
\min_{u_i \ge 0 \ v_i \ge 0}{\frac{1}{2}\|A - \sum_{i=1}^{r}{u_iv_i^T}\|^2_F} \label{NMF1}.
$$
\end{Problem}

Let us fix all the variables, except for a single vector $v_t$ and
consider the following least squares problem:
\begin{equation}
\min_{v \ge 0}{\frac{1}{2}\|R_t - u_tv^T\|^2_F} \label{NMF2},
\end{equation}

where $R_t = A - \sum_{i \neq t}{u_iv_i^T}$. We have:
\begin{eqnarray}
\|R_t - u_tv^T\|^2_F &=& trace\left[(R_t - u_tv^T)^T(R_t - u_tv^T)\right] \\
&=& \|R_t\|^2_F - 2v^TR_t^Tu_t + \|u_t\|^2_2\|v\|^2_2 \label{obj1}.
\end{eqnarray}
From this formulation, one now derives the following lemma.

\begin{Lemma} \label{ElementaryNMF}
If $[R_t^Tu_t]_+ \neq 0$, then $v_*:=\frac{[R_t^Tu_t]_+}{\|u_t\|_2^2}$ is the unique global minimizer of (\ref{NMF2}) and the function value equals $\|R_t\|^2_F-\frac{\|[R_t^Tu_t]_+\|^2_2}{\|u_t\|_2^2}.$
\end{Lemma}
\begin{proof}
Let us permute the elements of the vectors $x:=R_t^Tu_t$ and $v$ such that 
$$ Px= \left( \begin{array}{c} x_1 \\ x_2 \end{array}\right), \quad Pv= \left( \begin{array}{c} v_1 \\ v_2 \end{array}\right), \quad \mathrm{with} \quad x_1\ge 0, \quad x_2 < 0$$ 
and $P$ is the permutation matrix.
Then $$\|R_t - u_tv^T\|^2_F = \|R_t\|^2_F - 2v_1^Tx_1 - 2v_2^Tx_2  + \|u_t\|^2_2(v_1^Tv_1 + v_2^Tv_2).$$
Since $x_2<0$ and $v_2\geq 0$, it is obvious that $\|R_t - u_tv^T\|^2_F$ can only be minimal if $v_2=0$.
Our assumption implies that $x_1$ is nonempty and $x_1>0$. Moreover $[R_t^Tu_t]_+ \neq 0$ and $u_t \ge 0$ imply $\|u_t\|_2^2 > 0$, one can then find the optimal $v_1$ by minimizing the remaining quadratic function 
$$ \|R_t\|^2_F - 2v_1^Tx_1 + \|u_t\|^2_2v_1^Tv_1$$
which yields the solution $v_1=\frac{x_1}{\|u_t\|_2^2}$. Putting the two components together yields the result
$$ v_* = \frac{[R_t^Tu_t]_+}{\|u_t\|_2^2} \quad \text{and} \quad \|R_t - u_tv_*^T\|^2_F= \|R_t\|^2_F
-\frac{\|[R_t^Tu_t]_+\|^2_2}{\|u_t\|_2^2}.
$$
\end{proof}

\begin{algorithm}[htb]
\caption{(RRI)}
\label{RRI}
\begin{algorithmic}[1]
\STATE Initialize $u_i$'s, $v_i$'s, for $i=1$ to $r$
\REPEAT
\FOR {$t = 1$ to $r$}
\STATE $R_t=A - \sum_{i\neq t} {u_iv_i^T}$
\STATE{}
\IF{$[R_t^Tu_t]_+ \neq 0$}
	\STATE $v_t \leftarrow \frac{[R_t^Tu_t]_+}{\|u_t\|_2^2}$ \label{vupd}
\ELSE
	\STATE $v_t = 0$
\ENDIF
\STATE{}
\IF{$[R_tv_t]_+ \neq 0$}
	\STATE $u_t \leftarrow \frac{[R_tv_t]_+}{\|v_t\|_2^2}$ \label{uupd}
\ELSE
	\STATE $u_t = 0$
\ENDIF
\ENDFOR
\UNTIL{Stopping condition}
\end{algorithmic}
\end{algorithm}

\textbf{Remark 1}:
The above lemma has of course a dual form, where one fixes $v_t$ but solves for the optimal $u$
to minimize $\|R_t - uv_t^T\|^2_F$. This would yield the updating rules
\begin{equation} \label{updateui}
v_t \leftarrow \frac{[R_t^Tu_t]_+}{\|u_t\|_2^2} \quad \text{and} \quad
u_t \leftarrow \frac{[R_tv_t]_+}{\|v_t\|_2^2}
\end{equation}
which can be used to recursively update approximations $\sum_{i=1}^r u_iv_i^T$ by modifying each rank-one matrix $u_tv_t^T$ in a cyclic manner. This problem is different from the NMF, since the error matrices $R_t=A-\sum_{i\neq t} u_iv_i^T$ are no longer nonnegative. We will therefore call this method the \emph{Rank-one Residue Iteration} (RRI), i.e. Algorithm \ref{RRI}. The same algorithm was independently reported as Hierarchical Alternating Least Squares (HALS) \cite{Cichocki2007HALS}.

\textbf{Remark 2}: In case where $[R_t^Tu_t]_+ = 0$, we have a trivial solution for $v = 0$ that is not covered by Lemma \ref{ElementaryNMF}. In addition, if $u_t = 0$, this solution is no longer unique. In fact, $v$ can be arbitrarily taken to construct a rank-deficient approximation. The effect of this on the convergence of the algorithm will be discussed further in the next section. 

\textbf{Remark 3}: Notice that the optimality of Lemma \ref{ElementaryNMF} implies that $\|A - UV^T\|$ can not increase. And since $A \ge 0$ fixed,  $UV^T \ge 0$ must be bounded. Therefore, its component $u_iv_i^t$ (i=1\dots r) must be bounded as well. One can moreover scale the vector pairs $(u_i,v_i)$ at each stage as explained in Section \ref{scalingstopcond} without affecting the local optimality of Lemma \ref{ElementaryNMF}. It then follows that the rank one products $u_iv_i^T$ and their scaled vectors remain bounded.

\subsection{Convergence}

In the previous section, we have established the partial updates for each of the variable $u_i$ or $v_i$. And for a NMF problem where the reduced rank is $r$, we have in total $2r$ vector variables (the $u_i$'s and $v_i$'s). The described algorithm can be also considered as a projected gradient method since the update (\ref{updateui}) can be rewritten as:
\begin{eqnarray*}
u_t &\leftarrow& \frac{[R_tv_t]_+}{\|v_t\|_2^2} = \frac{[(A - \sum_{i\neq t}{u_iv_i^T})v_t]_+}{\|v_t\|_2^2} = \frac{[(A - \sum_{i}{u_iv_i^T} + u_tv_t^T)v_t]_+}{\|v_t\|_2^2}\\
&=& \frac{[(A - \sum_{i}{u_iv_i^T})v_t + u_tv_t^Tv_t]_+}{\|v_t\|_2^2} = \left[u_t - \frac{1}{\|v_t\|_2^2}\nabla_{u_t}\right]_+.
\end{eqnarray*}
Similarly, the update for $v_i$ can be rewritten as
$$
v_t \leftarrow \left[v_t - \frac{1}{\|u_t\|_2^2}\nabla_{v_t}\right]_+.
$$

Therefore, the new method follows the projected gradient scheme described in the previous section. But it  produces the optimal solution in closed form. For each update of a column $v_t$ (or $u_t$), the proposed algorithm requires just a matrix-vector multiplication $R_t^Tu_t$ (or $R_tv_t$), wherein the residue matrix $R_t=A - \sum_{i\neq t} {u_iv_i^T}$ does not have to be calculated explicitly. Indeed, by calculating $R_t^Tu_t$ (or $R_tv_t$) from $A^Tu_t$ (or $Av_t$) and $\sum_{i\neq t} {v_i(u_i^T u_t)}$ (or $\sum_{i\neq t} {u_i(v_i^T v_t)}$), the complexity is reduced from $O(mnr + mn)$ to only $O\left (mn + (m+n)(r-1) \right )$ which is majored by $O(mn)$. This implies that the complexity of each sweep through the $2r$ variables $u_t's$ and $v_t's$ requires only $O(mnr)$ operations, which is equivalent to a sweep of the multiplicative rules and to an inner loop of any gradient methods. This is very low since the evaluation of the whole gradient requires already the same complexity. 

Because at each step of the 2$r$ basic steps of Algorithm \ref{RRI}, we compute an optimal rank-one nonnegative correction to the corresponding error matrix $R_t$ the Frobenius norm of the error can not increase. This is a reassuring property but it does not imply convergence of the algorithm.

Each vector $u_t$ or $v_t$ lies in a convex set $\mathbb{U}_t \subset \mathbb{R}_+^m$ or $\mathbb{V}_t \subset \mathbb{R}_+^n$. Moreover, because of the possibility to include scaling we can set an upper bound for $\|U\|$ and $\|V\|$, in such a way that all the $\mathbb{U}_t$ and $\mathbb{V}_t$ sets can be considered as closed convex. Then, we can use the following Theorem \ref{convergence}, to prove a stronger convergence result for Algorithm \ref{RRI}.

\begin{Theorem} \label{convergence}
Every limit point generated by Algorithm \ref{RRI} is a stationary point.
\end{Theorem}
\begin{proof}
We notice that, if $u_t = 0$ and $v_t = 0$ at some stages of Algorithm \ref{RRI}, they will remain zero and no longer take part in all subsequent iterations. We can divide the execution of Algorithm \ref{RRI} into two phases. 

During the first phase, some of the pairs $(u_t, v_t)$ become zero. Because there are only a finite number ($2r$) of such vectors, the number of iterations in this phase is also finite. At the end of this phase, we can rearrange and partition the matrices $U$ and $V$ such that
$$
U = (U_+ \ 0) \ \text{and} \ V = (V_+ \ 0),
$$
where $U_+$ and $V_+$ do not have any zero column. We temporarily remove zero columns out of the approximation.

During the second phase, no column of $U_+$ and $V_+$ becomes zero, which guarantees the updates for the columns of $U_+$ and $V_+$ are unique and optimal. Moreover, $\frac{1}{2}\|A - \sum_{i=1}^{r}{u_iv_i^T}\|^2_F$ is continuously differentiable over the set $\mathbb{U}_1\times \ldots \times \mathbb{U}_r \times \mathbb{V}_1 \times \ldots \times \mathbb{V}_r$, and the $\mathbb{U}_i$'s and $\mathbb{V}_i$'s are closed convex. A direct application of Proposition 2.7.1 in \cite{bertsekas1999np} proves that every stationary point $(U^*_+, V^*_+)$ is a stationary point. It is then easy to prove that if there are zero columns removed at the end of the first phase, adding them back yields another stationary point: $U^* = (U^*_+ \ 0)$ and $V^* = (V^*_+ \ 0)$ of the required dimension. However, in this case, the rank of the approximation will then be lower than the requested dimension $r$.
\end{proof}

In Algorithm \ref{RRI}, variables are updated in this order: $u_1$, $v_1$, $u_2$, $v_2$, $\ldots$. We can alternate the variables in a different order as well, for example $u_1$, $u_2$, $\ldots$, $u_r$ $v_1$, $v_2$, $\ldots$, $v_r$, $\ldots$. Whenever this is carried out in a cyclic fashion, the Theorem \ref{convergence} still holds and this does not increase the complexity of each iteration of the algorithm.

As pointed above, stationary points given by Algorithm \ref{RRI} may contain useless zero components. To improve this, one could replace $u_tv_t^T (\equiv 0)$ by any nonnegative rank-one approximation that reduces the norm of the error matrix. For example, the substitution
\begin{equation}
u_t = e_{i^*} \qquad v_t = [R_t^Tu_t]_+,  \label{r1substitute}
\end{equation}
where $i^* = \argmax_{i}{\|[R_t^Te_i]_+\|_2^2}$, reduces the error norm by $\|[R_t^Te_i]_+\|_2^2 > 0$ unless $R_t \le 0$. These substitutions can be done as soon as  $u_t$ and $v_t$ start to be zero. If we do these substitutions in only a \emph{finite} number of times before the algorithm starts to converge, Theorem \ref{convergence} still holds. In practice, only a few such substitutions in total are usually needed by the algorithm to converge to a stationary point without any zero component. Note that the matrix rank of the approximation might not be $r$, even when all $u_t$'s and $v_t$'s ($t = 1\dots r$) are nonzero. 

A possibly better way to fix the problem due to zero components is to use the following \emph{damped RRI algorithm} in which we introduce new $2r$ dummy variables $w_i \in \mathbb U_i$ and $z_i \in \mathbb V_i$, where $i = 1...r$. The new problem to solve is:
\begin{Problem} [Damped Nonnegative Matrix Factorization] \label{dNMF}
$$
\min_{\substack{u_i \ge 0 \ v_i \ge 0 \\ w_i \ge 0 \ z_i \ge 0}}{\frac{1}{2}\|A - \sum_{i=1}^{r}{u_iv_i^T}\|^2_F + \frac{\psi}{2}\sum_i \|u_i - w_i\|^2_2 + \frac{\psi}{2}\sum_i \|v_i - z_i\|^2_2},
$$
where the damping factor $\psi$ is a positive constant.
\end{Problem}

Again, the coordinate descent scheme is applied with the cyclic update order: $u_1$, $w_1$, $v_1$, $z_1$, $u_2$, $w_2$, $v_2$, $z_2$, $\ldots$ to result in the following optimal updates for $u_t$, $v_t$, $w_t$ and $z_t$:
\begin{equation}
u_t = \frac{[R_tv_t]_+ + \psi w_t}{\|v_t\|_2^2 + \psi}, \ w_t = u_t, \ v_t = \frac{[R_t^Tu_t]_+ + \psi z_t}{\|u_t\|_2^2 + \psi} \ \text{and} \ z_t = v_t \label{drri}
\end{equation}
where $t = 1\dots r$. The updates $w_t = u_t$ and $z_t = v_t$ can be integrated in the updates of $u_t$ and $v_t$ to yield Algorithm \ref{alg:dampedRRI}. We have the following results:
\begin{Theorem} \label{convergencedrri}
Every limit point generated by Algorithm \ref{alg:dampedRRI} is a stationary point of NMF problem \ref{NMF1}.
\end{Theorem}

\begin{algorithm}[tph]
\caption{(Damped RRI)}
\label{alg:dampedRRI}
\begin{algorithmic}[1]
\STATE Initialize $u_i$'s, $v_i$'s, for $i=1$ to $r$
\REPEAT
\FOR {$t = 1$ to $r$}
\STATE $R_t=A - \sum_{i\neq t} {u_iv_i^T}$
\STATE $v_t \leftarrow \frac{[R_t^Tu_t + \psi v_t]_+}{\|u_t\|_2^2 + \psi}$ 
\STATE $u_t \leftarrow \frac{[R_tv_t + \psi u_t]_+}{\|v_t\|_2^2 + \psi}$ 
\ENDFOR
\UNTIL{Stopping condition}
\end{algorithmic}
\end{algorithm}

\begin{proof}
Clearly the cost function in Problem \ref{dNMF} is continuously differentiable over the set $\mathbb{U}_1\times \ldots \times \mathbb{U}_r \times \mathbb{U}_1\times \ldots \times \mathbb{U}_r$ $\times \mathbb{V}_1 \times \ldots \times \mathbb{V}_r \times \mathbb{V}_1 \times \ldots \times \mathbb{V}_r$, and the $\mathbb{U}_i$'s and $\mathbb{V}_i$'s are closed convex. The uniqueness of the global minimum of the elementary problems and a direct application of Proposition 2.7.1 in \cite{bertsekas1999np} prove that every limit point of Algorithm \ref{alg:dampedRRI} is a stationary point of Problem \ref{dNMF}.

Moreover, at a stationary point of Problem \ref{dNMF}, we have $u_t = w_t$ and $v_t = z_t$, $t=1...r$. The cost function in Problem \ref{dNMF} becomes the cost function of the NMF problem \ref{NMF1}. This implies that every stationary point of Problem \ref{dNMF} yields a stationary point of the standard NMF problem \ref{NMF1}.
\end{proof}

This \emph{damped} version not only helps to eliminate the problem of zero components in the convergence analysis but may also help to avoid zero columns in the approximation when $\psi$ is carefully chosen. But it is not an easy task. Small values of $\psi$ provide an automatic treatment of zeros while not changing much the updates of RRI. Larger values of $\psi$ might help to prevent the vectors $u_t$ and $v_t$ ($t = 1 \dots r$) from becoming zero too soon. But too large values of $\psi$ limit the updates to only small changes, which will slow down the convergence. 

In general, the rank of the approximation can still be lower than the requested dimension. Patches may still be needed when a zero component appears. Therefore, in our experiments, using the \emph{undamped} RRI algorithm \ref{RRI} with the substitution (\ref{r1substitute}) is still the best choice. 

\subsection{Variants of the RRI method} \label{RRIVars}

We now extend the Rank-one Residue Iteration by using a factorization of the type $XDY^T$ where $D$ is diagonal and nonnegative and the columns of the nonnegative matrices $X$ and $Y$ are normalized. The NMF formulation then becomes
$$
\min_{\substack{x_i \in \mathbb{X}_i \ y_i \in \mathbb{Y}_i \\ d_i \in \mathbb{R}_+}}{\frac{1}{2}\|A - \sum_{i=1}^{r}{d_ix_iy_i^T}\|^2_F} \label{NMFVar},
$$
where $\mathbb{X}_i$'s and  $\mathbb{Y}_i$'s are sets of normed vectors. 

The variants that we present here depend on the choice of $\mathbb{X}_i$'s and  $\mathbb{Y}_i$'s. 
A generalized Rank-one Residue Iteration method for low-rank approximation is given in Algorithm \ref{GRRI}. This algorithm needs to solve a sequence of elementary problems of the type:
\begin{equation}
\max_{s \in \mathbb{S}}{y^Ts} \label{maxinprod}
\end{equation}
where $y \in \mathbb{R}^n$ and $\mathbb{S} \subset \mathbb{R}^n$ is a set of normed vectors. We 
first introduce a permutation vector $I_y = (i_1 \ i_2 \ \ldots \ i_n)$ which reorders the elements of
$y$ in non-increasing order : $y_{i_k} \ge y_{i_{k+1}}$, $k = 1\ldots(n-1)$. The function $p(y)$ returns the number of positive entries of $y$.

\begin{algorithm}[tph]
\caption{GRRI}
\label{GRRI}
\begin{algorithmic}[1]
\STATE Initialize $x_i$'s, $y_i$'s and $d_i$'s, for $i=1$ to $r$
\REPEAT
\FOR {$i = 1$ to $r$}
\STATE $R_i=A - \sum_{j\neq i} {d_jx_jy_j^T}$
\STATE $y_i \leftarrow \argmax_{s \in \mathbb{Y}_i}{\left( x_i^TR_is \right)}$
\STATE $x_i \leftarrow \argmax_{s \in \mathbb{X}_i}{\left( y_i^TR_i^Tc \right)}$
\STATE $d_i = x_i^TR_iy_i$
\ENDFOR
\UNTIL{Stopping condition}
\end{algorithmic}
\end{algorithm}

Let us first point out that for the set of normed nonnegative vectors the solution of 
problem (\ref{maxinprod}) is given by $s^* = \frac{y_+}{\|y_+\|_2}$. It then follows that
Algorithm \ref{GRRI} is essentially the same as Algorithm \ref{RRI} since the solutions 
$v_i$ and $u_i$ of each step of Algorithm \ref{GRRI}, given by (\ref{updateui}), 
correspond exactly to those of problem (\ref{maxinprod}) via the relations
$y_i=u_i/\|u_i\|_2$, $y_i=v_i/\|v_i\|_2$ and $d_i=\|u_i\|_2\|v_i\|_2$.

Below we list the sets for which the solution $s^*$ of (\ref{maxinprod}) can be easily computed.
\begin{itemize}
\item \textit{Set of normed vectors}: $s = \frac{y}{\|y\|_2}$. This is useful when one wants to create factorizations where only one of the factor $U$ or $V$ is nonnegative and the other is real matrix.
\item \textit{Set of normed nonnegative vectors}: $s = \frac{y_+}{\|y_+\|_2}$.  
\item \textit{Set of normed bounded nonnegative vectors $\{s\}$}: where $0 \le l_i \le s_i \le p_i$. The optimal solution of (\ref{maxinprod}) is given by:
$$
s = \max \left( l, \ \min \left ( p, \ \frac{y_+}{\|y_+\|_2} \right ) \right).
$$
\item \textit{Set of normed binary vectors $\{s\}$}: where $s = \frac{b}{\|b\|}$ and $b \in \{0, 1\}^n$.  The optimal solution of (\ref{maxinprod}) is given by:
$$
[s^*]_{i_t} = \left \{ \begin{array}{ll} \frac{1}{\sqrt{k^*}} & \text{if } t \le k^* \\ 0 & \text{otherwise} \end{array} \right. \quad \text{where} \
k^* = \argmax_{k}{\frac{\sum_{t=1}^{k}y_{i_t}}{\sqrt k}}.
$$
\item \textit{Set of normed sparse nonnegative vectors}: all normed nonnegative vectors having at most $K$ nonzero entries. The optimal solution for (\ref{maxinprod}) is given by norming the following vector $p^*$
$$
[p^*]_{i_t} = \left \{ \begin{array}{ll} y_{i_t} & \text{if } t \le \min(p(y), K) \\ 0 & \text{otherwise} \end{array} \right.
$$
\item \textit{Set of normed fixed-sparsity nonnegative vectors}: all nonnegative vectors $s$ a fixed sparsity, where
$$
sparsity(s) = \frac{\sqrt{n} - \|s\|_1/\|s\|_2}{\sqrt{n} - 1}.
$$
The optimal solution for (\ref{maxinprod}) is given by using the projection scheme in \cite{hoyer2004}.
\end{itemize}

One can also imagine other variants, for instance by combining the above ones. Depending on how data 
need to be approximated, one can create new algorithms provided it is relatively simple to solve problem (\ref{maxinprod}). There have been some particular ideas in the literatures such as NMF with sparseness constraint \cite{hoyer2004}, Semidiscrete Matrix Decomposition \cite{kolda1998smd} and Semi-Nonnegative Matrix Factorization \cite{ding2006seminmf} for which variants of the above scheme can offer an alternative choice of algorithm.

\textbf{Remark:} Only the first three sets are the normed version of a closed convex set, as required for the convergence by Theorem \ref{convergence}. Therefore the algorithms might not converge to a stationary point with the other sets. However, the algorithm always guarantees a non-increasing update even in those cases and can therefore be expected to return a \emph{good} approximation.

\subsection{Nonnegative Tensor Factorization}

If we refer to the problem of finding the nearest nonnegative vector to a given vector $a$ as the nonnegative approximation in one dimension, the NMF is its generalization in two dimensions and naturally, it can be extended to even higher-order tensor approximation problems. Algorithms described in the previous section use the closed form solution of the one dimensional problem to solve the two-dimensional problem. We now generalize this to higher orders. Since in one dimension such an approximation is easy to construct, we continue to use this approach to build the solutions for higher order problems.

For a low-rank tensor, there are two popular kinds of factored tensors, namely those of Tucker and Kruskal \cite{badertensor}. We only give an algorithm for finding approximations of Kruskal type. It is easy to extend this to tensors of Tucker type, but this is omitted here.

Given a $d$ dimensional tensor $T$, we will derive an algorithm for approximating a nonnegative tensor by a rank-$r$ nonnegative Kruskal tensor $S \in \mathbb{R}_+^{n_1 \times n_2 \times \ldots \times n_d}$ represented as a sum of $r$ rank-one tensors:
$$
S = \sum_{i = 1}^{r}{\sigma_i u_{1i} \star u_{2i} \star \ldots \star u_{di}}
$$
where $\sigma_i \in \mathbb{R}_+$ is a scaling factor, $u_{ti} \in \mathbb{R}_+^{n_t}$ is a normed vector (i.e. $\|u_{ti}\|_2 = 1$) and $a \star b$ stands for the outer product between two vectors or tensors $a$ and $b$. 

The following update rules are the generalization of the matrix case to the higher order tensor:
\begin{eqnarray}
     y & = & (\ldots((\ldots(R_k u_{1k})\ldots u_{(t-1)k})u_{(t+1)k}) \ldots )u_{dk}  \\
    \sigma_k & = & \| [y]_+ \|_2 , \quad u_{tk} = \frac{[y]_+}{\sigma_k},
\end{eqnarray}
where $R_k = T - \sum_{i \neq k}{\sigma_i u_{1i} \star u_{2i} \star \ldots \star u_{di}}$ is the residue tensor calculated without the $k^{th}$ component of $S$ and $R_ku_{ij}$ is the ordinary tensor/vector product in the corresponding dimension.

We can then produce an algorithm which updates in a cyclic fashion all vectors $u_{ji}$. This is in fact a direct extension to Algorithm \ref{RRI}, one can carry out the same discussion about the convergence here to guarantee that each limit point of this algorithm is a stationary point for the nonnegative tensor factorization problem and to improve the approximation quality.

Again, as we have seen in the previous section, we can extend the procedure to take into account different constraints on the vectors $u_{ij}$ such as discreteness, sparseness, etc.

The approach proposed here is again different from that in \cite{cichocki2008nma} where a similar cascade procedure for multilayer nonnegative matrix factorization is used to compute a 3D tensor approximation. Clearly, the approximation error will be higher than our proposed method, since the cost function is minimized by taking into account all the dimensions.

\subsection{Regularizations} \label{regularize}

The regularizations are common methods to cope with the ill-posedness of inverse problems. Having known some additional information about the solution, one may want to imposed a priori some constraints to algorithms, such as: smoothness, sparsity, discreteness, etc. To add such regularizations in to the RRI algorithms, it is possible to modify the NMF cost function by adding some regularizing terms. We will list here the update for $u_i$'s and $v_i$'s when some simple regularizations are added to the original cost function. The proof of these updates are straight-forward and hence omitted.  

\begin{itemize}
\item One-Norm $\|.\|_1$ regularization: the one-norm of the vector variable can be added as a heuristic for finding a sparse solution. This is an alternative to the fixed-sparsity variant presented above. The regularized cost function with respect to the variable $v_t$ will be 
$$
\frac{1}{2}\|R_t - u_tv^T\|^2_F + \beta \|v\|_1, \quad \beta > 0
$$
where the optimal update is given by
$$
v^*_t = \frac{[R_t^Tu_t - \beta \mathbf{1}_{n \times 1}]_+}{\|u_t\|_2^2}.
$$
The constant $\beta > 0$ can be varied to control the trade-off between the approximation error $\frac{1}{2}\|R_t - u_tv^T\|^2_F$ and $\|v\|_1$. From this update, one can see that this works by zeroing out elements of $R_t^Tu_t$ which are smaller than $\beta$, hence reducing the number of nonzero elements of $v^*_t$.

\item Smoothness regularization $\|v - B\hat v_t\|_F^2$: where $\hat v_t$ is the current value of $v_t$ and the matrix $B$ helps to calculate the average of the neighboring elements at each element of $v$. When $v$ is a 1D smooth function, $B$ can be the following $n \times n$ matrix:
\begin{equation}
B = \left (\begin{array}{ccccc} 
0 & 1 & \ldots & \ldots & 0  \\
\frac{1}{2} & 0 & \frac{1}{2} & \ldots & 0  \\
\vdots & \ddots & \ddots & \ddots & \vdots  \\
0  & \ldots  & \frac{1}{2} & 0 & \frac{1}{2} \\
0  & \ldots & 0 & 1 & 0 
\end{array} \right ). \label{smoothmatrix}
\end{equation}
This matrix can be defined in a different way to take the true topology of $v$ into account, for instance $v = vec(F)$ where $F$ is a matrix. The regularized cost function with respect to the variable $v_t$ will be 
$$
\frac{1}{2}\|R_t - u_tv^T\|^2_F + \frac{\delta}{2} \|v - B\hat v_t\|_F^2, \quad \delta > 0
$$
where the optimal update is given by
$$
v^*_t = \frac{[R_t^Tu_t + \delta B\hat v_t]_+}{\|u_t\|_2^2 + \delta}.
$$
The constant $\delta \ge 0$ can be varied to control the trade-off between the approximation error $\frac{1}{2}\|R_t - u_tv^T\|^2_F$ and the smoothness of $v_t$ at the fixed point. From the update, one can see that this works by searching for the optimal update $v^*_t$ with some preference for the neighborhood of $B\hat v_i$, i.e., a smoothed vector of the current value $\hat v_t$.
\end{itemize}

The two above regularizations can be added independently to each of the columns of $U$ and/or $V$. The trade-off factor $\beta$ (or $\delta$) can be different for each column. A combination of different regularizations on a column (for instance $v_t$) can also be used to solve the multi-criterion problem
$$
\frac{1}{2}\|R_t - u_tv^T\|^2_F  + \frac{\gamma}{2} \|v\|_2^2 + \frac{\delta}{2} \|v - B\hat v_t\|_F^2, \quad \beta, \gamma, \delta > 0
$$
where the optimal update is given by 
$$
v^*_t = \frac{[R_t^Tu_t - \beta \mathbf{1}_{n \times 1} + \delta B\hat v_t]_+}{\|u_t\|_2^2 + \delta}.
$$

The one-norm regularizations as well as the two-norm regularization can be found in \cite{albright2006aia} and \cite{berry2007aaa}. A major difference with that method is that the norm constraints is added to the rows rather than on the columns of $V$ or $U$ as done here. However, for the two versions of the one-norm regularization, the effects are somehow similar. While the two-norm regularization on the columns of $U$ and $V$ are simply scaling effects, which yield nothing in the RRI algorithm. We therefore only test the smoothness regularization at the end of the chapter with some numerical generated data.

For more extensions and variants, see \cite{ho2008thesis}.

\section{Experiments} \label{experiment}

Here we present several experiments to compare the different descent algorithms presented in this paper. For all the algorithms, the scaling scheme proposed in section \ref{scalingstopcond} was applied.

\subsection{Random matrices}

\begin{table}[tbp]
\begin{center}
\begin{tabular}{crrrrrrr}
\hline
 $\epsilon$& Mult& ALS& FLine& CLine& FFO& CFO& RRI\\ 
\hline 
\multicolumn{8}{c}{$\mathbf{{\scriptstyle (m = 30, \ n = 20, \ r = 2)}}$} \\ 
${\scriptstyle 10^{-2}}$& ${\scriptstyle 0.02(96)}$& ${\scriptstyle 0.40}$& ${\scriptstyle 0.04}$& ${\scriptstyle 0.02}$& ${\scriptstyle 0.02}$& ${\scriptstyle 0.01}$& ${\scriptstyle 0.01}$\\ 
${\scriptstyle 10^{-3}}$& ${\scriptstyle 0.08(74)}$& ${\scriptstyle 1.36}$& ${\scriptstyle 0.12}$& ${\scriptstyle 0.09}$& ${\scriptstyle 0.05}$& ${\scriptstyle 0.04}$& ${\scriptstyle 0.03}$\\ 
${\scriptstyle 10^{-4}}$& ${\scriptstyle 0.17(71)}$& ${\scriptstyle 2.81}$& ${\scriptstyle 0.24}$& ${\scriptstyle 0.17}$& ${\scriptstyle 0.11}$& ${\scriptstyle 0.08}$& ${\scriptstyle 0.05}$\\ 
${\scriptstyle 10^{-5}}$& ${\scriptstyle 0.36(64)}$& ${\scriptstyle 4.10}$& ${\scriptstyle 0.31}$& ${\scriptstyle 0.25}$& ${\scriptstyle 0.15}$& ${\scriptstyle 0.11}$& ${\scriptstyle 0.07}$\\ 
${\scriptstyle 10^{-6}}$& ${\scriptstyle 0.31(76)}$& ${\scriptstyle 4.74}$& ${\scriptstyle 0.40}$& ${\scriptstyle 0.29}$& ${\scriptstyle 0.19}$& ${\scriptstyle 0.15}$& ${\scriptstyle 0.09}$\\ 
\hline 
\multicolumn{8}{c}{$\mathbf{{\scriptstyle (m = 100, \ n = 50, \ r = 5)}}$} \\ 
 ${\scriptstyle 10^{-2}}$& ${\scriptstyle 45*(0)}$& ${\scriptstyle 3.48}$& ${\scriptstyle 0.10}$& ${\scriptstyle 0.09}$& ${\scriptstyle 0.09}$& ${\scriptstyle 0.04}$& ${\scriptstyle 0.02}$\\ 
${\scriptstyle 10^{-3}}$& ${\scriptstyle 45*(0)}$& ${\scriptstyle 24.30(96)}$& ${\scriptstyle 0.59}$& ${\scriptstyle 0.63}$& ${\scriptstyle 0.78}$& ${\scriptstyle 0.25}$& ${\scriptstyle 0.15}$\\ 
${\scriptstyle 10^{-4}}$& ${\scriptstyle 45*(0)}$& ${\scriptstyle 45*(0)}$& ${\scriptstyle 2.74}$& ${\scriptstyle 2.18}$& ${\scriptstyle 3.34}$& ${\scriptstyle 0.86}$& ${\scriptstyle 0.45}$\\ 
${\scriptstyle 10^{-5}}$& ${\scriptstyle 45*(0)}$& ${\scriptstyle 45*(0)}$& ${\scriptstyle 5.93}$& ${\scriptstyle 4.06}$& ${\scriptstyle 6.71}$& ${\scriptstyle 1.58}$& ${\scriptstyle 0.89}$\\ 
${\scriptstyle 10^{-6}}$& ${\scriptstyle 45*(0)}$& ${\scriptstyle 45*(0)}$& ${\scriptstyle 7.23}$& ${\scriptstyle 4.75}$& ${\scriptstyle 8.98}$& ${\scriptstyle 1.93}$& ${\scriptstyle 1.30}$\\ 
\hline 
\multicolumn{8}{c}{$\mathbf{{\scriptstyle (m = 100, \ n = 50, \ r = 10)}}$} \\ 
${\scriptstyle 10^{-2}}$& ${\scriptstyle 45*(0)}$& ${\scriptstyle 11.61}$& ${\scriptstyle 0.28}$& ${\scriptstyle 0.27}$& ${\scriptstyle 0.18}$& ${\scriptstyle 0.11}$& ${\scriptstyle 0.05}$\\ 
${\scriptstyle 10^{-3}}$& ${\scriptstyle 45*(0)}$& ${\scriptstyle 41.89(5)}$& ${\scriptstyle 1.90}$& ${\scriptstyle 2.11}$& ${\scriptstyle 1.50}$& ${\scriptstyle 0.74}$& ${\scriptstyle 0.35}$\\ 
${\scriptstyle 10^{-4}}$& ${\scriptstyle 45*(0)}$& ${\scriptstyle 45*(0)}$& ${\scriptstyle 7.20}$& ${\scriptstyle 5.57}$& ${\scriptstyle 5.08}$& ${\scriptstyle 2.29}$& ${\scriptstyle 1.13}$\\ 
${\scriptstyle 10^{-5}}$& ${\scriptstyle 45*(0)}$& ${\scriptstyle 45*(0)}$& ${\scriptstyle 12.90}$& ${\scriptstyle 9.69}$& ${\scriptstyle 10.30}$& ${\scriptstyle 4.01}$& ${\scriptstyle 1.71}$\\ 
${\scriptstyle 10^{-6}}$& ${\scriptstyle 45*(0)}$& ${\scriptstyle 45*(0)}$& ${\scriptstyle 14.62(99)}$& ${\scriptstyle 11.68(99)}$& ${\scriptstyle 13.19}$& ${\scriptstyle 5.26}$& ${\scriptstyle 2.11}$\\ 
\hline 
\multicolumn{8}{c}{$\mathbf{{\scriptstyle (m = 100, \ n = 50, \ r = 15)}}$} \\ 
${\scriptstyle 10^{-2}}$& ${\scriptstyle 45*(0)}$& ${\scriptstyle 25.98}$& ${\scriptstyle 0.66}$& ${\scriptstyle 0.59}$& ${\scriptstyle 0.40}$& ${\scriptstyle 0.20}$& ${\scriptstyle 0.09}$\\ 
${\scriptstyle 10^{-3}}$& ${\scriptstyle 45*(0)}$& ${\scriptstyle 45*(0)}$& ${\scriptstyle 3.90}$& ${\scriptstyle 4.58}$& ${\scriptstyle 3.18}$& ${\scriptstyle 1.57}$& ${\scriptstyle 0.61}$\\ 
${\scriptstyle 10^{-4}}$& ${\scriptstyle 45*(0)}$& ${\scriptstyle 45*(0)}$& ${\scriptstyle 16.55(98)}$& ${\scriptstyle 13.61(99)}$& ${\scriptstyle 9.74}$& ${\scriptstyle 6.12}$& ${\scriptstyle 1.87}$\\ 
${\scriptstyle 10^{-5}}$& ${\scriptstyle 45*(0)}$& ${\scriptstyle 45*(0)}$& ${\scriptstyle 21.72(97)}$& ${\scriptstyle 17.31(92)}$& ${\scriptstyle 16.59(98)}$& ${\scriptstyle 7.08}$& ${\scriptstyle 2.39}$\\ 
${\scriptstyle 10^{-6}}$& ${\scriptstyle 45*(0)}$& ${\scriptstyle 45*(0)}$& ${\scriptstyle 25.88(89)}$& ${\scriptstyle 19.76(98)}$& ${\scriptstyle 19.20(98)}$& ${\scriptstyle 10.34}$& ${\scriptstyle 3.66}$\\ 
\hline 
\multicolumn{8}{c}{$\mathbf{{\scriptstyle (m = 100, \ n = 100, \ r = 20)}}$} \\  
${\scriptstyle 10^{-2}}$& ${\scriptstyle 45*(0)}$& ${\scriptstyle 42.51(4)}$& ${\scriptstyle 1.16}$& ${\scriptstyle 0.80}$& ${\scriptstyle 0.89}$& ${\scriptstyle 0.55}$& ${\scriptstyle 0.17}$\\ 
${\scriptstyle 10^{-3}}$& ${\scriptstyle 45*(0)}$& ${\scriptstyle 45*(0)}$& ${\scriptstyle 9.19}$& ${\scriptstyle 8.58}$& ${\scriptstyle 10.51}$& ${\scriptstyle 5.45}$& ${\scriptstyle 1.41}$\\ 
${\scriptstyle 10^{-4}}$& ${\scriptstyle 45*(0)}$& ${\scriptstyle 45*(0)}$& ${\scriptstyle 28.59(86)}$& ${\scriptstyle 20.63(94)}$& ${\scriptstyle 29.89(69)}$& ${\scriptstyle 12.59}$& ${\scriptstyle 4.02}$\\ 
${\scriptstyle 10^{-5}}$& ${\scriptstyle 45*(0)}$& ${\scriptstyle 45*(0)}$& ${\scriptstyle 32.89(42)}$& ${\scriptstyle 27.94(68)}$& ${\scriptstyle 34.59(34)}$& ${\scriptstyle 18.83(90)}$& ${\scriptstyle 6.59}$\\ 
${\scriptstyle 10^{-6}}$& ${\scriptstyle 45*(0)}$& ${\scriptstyle 45*(0)}$& ${\scriptstyle 37.14(20)}$& ${\scriptstyle 30.75(60)}$& ${\scriptstyle 36.48(8)}$& ${\scriptstyle 22.80(87)}$& ${\scriptstyle 8.71}$\\ 
\hline 
\multicolumn{8}{c}{$\mathbf{{\scriptstyle (m = 200, \ n = 100, \ r = 30)}}$} \\ 
${\scriptstyle 10^{-2}}$& ${\scriptstyle 45*(0)}$& ${\scriptstyle 45*(0)}$& ${\scriptstyle 2.56}$& ${\scriptstyle 2.20}$& ${\scriptstyle 2.68}$& ${\scriptstyle 1.31}$& ${\scriptstyle 0.44}$\\ 
${\scriptstyle 10^{-3}}$& ${\scriptstyle 45*(0)}$& ${\scriptstyle 45*(0)}$& ${\scriptstyle 22.60(99)}$& ${\scriptstyle 25.03(98)}$& ${\scriptstyle 29.67(90)}$& ${\scriptstyle 12.94}$& ${\scriptstyle 4.12}$\\ 
${\scriptstyle 10^{-4}}$& ${\scriptstyle 45*(0)}$& ${\scriptstyle 45*(0)}$& ${\scriptstyle 36.49(2)}$& ${\scriptstyle 39.13(13)}$& ${\scriptstyle 45*(0)}$& ${\scriptstyle 33.33(45)}$& ${\scriptstyle 14.03}$\\ 
${\scriptstyle 10^{-5}}$& ${\scriptstyle 45*(0)}$& ${\scriptstyle 45*(0)}$& ${\scriptstyle 45*(0)}$& ${\scriptstyle 39.84(2)}$& ${\scriptstyle 45*(0)}$& ${\scriptstyle 37.60(6)}$& ${\scriptstyle 21.96(92)}$\\ 
${\scriptstyle 10^{-6}}$& ${\scriptstyle 45*(0)}$& ${\scriptstyle 45*(0)}$& ${\scriptstyle 45*(0)}$& ${\scriptstyle 45*(0)}$& ${\scriptstyle 45*(0)}$& ${\scriptstyle 45*(0)}$& ${\scriptstyle 25.61(87)}$\\ 
\hline 
\end{tabular}

\caption{Comparison of average successful running time of algorithms over $100$ random matrices. Time limit is $45$ seconds. $0.02(96)$ means that a result is returned with the required precision $\epsilon$ within $45$ seconds for $96$ (of $100$) matrices of which the average running time is $0.02$ seconds. $45*(0)$: failed in all $100$ matrices.} \label{tab:results}
 \end{center} 
 \end{table} 

We generated 100 random nonnegative matrices of different sizes. We used seven different algorithms to approximate each matrix:
\begin{itemize}
\item the multiplicative rule (\textbf{Mult}),
\item alternative least squares using Matlab function \textit{lsqnonneg} (\textbf{ALS}),
\item a full space search using line search and Armijo criterion (\textbf{FLine}),
\item a coordinate search alternating on $U$ and $V$, and using line search and Armijo criterion (\textbf{CLine}),
\item a full space search using first-order approximation (\textbf{FFO}),
\item a coordinate search alternating on $U$ and $V$, and using first-order approximation (\textbf{CFO})
\item an iterative rank-one residue approximation (\textbf{RRI}).
\end{itemize}

For each matrix, the same starting point is used for every algorithm. We create a starting point by randomly generating two matrices $U$ and $V$ and then rescaling them to yield a first approximation of the original matrix $A$ as proposed in Section \ref{scalingstopcond}:
$$
\qquad U = UD\sqrt{\alpha}, \qquad V = VD^{-1}\sqrt{\alpha},
$$
where 
$$
\alpha := \frac{\left<A, UV^T\right>}{\left<UV^T, UV^T\right>} \quad \text{and} \quad D_{ij} = \left \{ \begin{array}{cl} \sqrt{\frac{\|V_{:i}\|_2}{\|U_{:i}\|_2}} & \quad \text{if } i = j \\ 0 & \quad \text{otherwise} \end{array} \right. .
$$

From (\ref{sumkkt}), we see that when approaching a KKT stationary point of the problem, the above scaling factor $\alpha \rightarrow 1$. This implies that every KKT stationary point of this problem is scale-invariant.  

The algorithms are all stopped when the projected gradient norm is lower than $\epsilon$ times the gradient norm at the starting point or when it takes more than $45$ seconds. The relative precisions $\epsilon$ are chosen equal to $10^{-2}$, $10^{-3}$, $10^{-4}$, $10^{-5}$, $10^{-6}$. No limit was imposed on the number of iterations.

For alternative gradient algorithms \textbf{CLine} and \textbf{CFO}, we use different precisions $\epsilon_U$ and $\epsilon_V$ for each of the inner iteration for $U$ and for $V$ as suggested in \cite{cjl05b} where $\epsilon_U$ and $\epsilon_V$ are initialized by $10^{-3}$. And when the inner loop for $U$ or $V$ needs no iteration to reach the precision $\epsilon_U$ or $\epsilon_V$, one more digit of precision will be added into $\epsilon_U$ or $\epsilon_V$ (i.e. $\epsilon_U = \epsilon_U /10$ or $\epsilon_V = \epsilon_V /10$).

Table \ref{tab:results} shows that for all sizes and ranks, Algorithm \textbf{RRI} is the fastest to reach the required precision. Even though it is widely used in practice, algorithm \textbf{Mult} fails to provide solutions to the NMF problem within the allocated time. A further investigation shows that the algorithm gets easily trapped in boundary points where some $U_{ij}$ and/or $V_{ij}$ is zero while $\nabla_{U_{ij}}$ and/or $\nabla_{V_{ij}}$ is negative, hence violating one of the KKT conditions (\ref{KKT2}). The multiplicative rules then fail to move and do not return to a local minimizer. A slightly modified version of this algorithm was given in \cite{cjl05c}, but it needs to wait to get sufficiently close to such points before attempting an escape, and is therefore also not efficient. The \textbf{ALS} algorithm can return a stationary point, but it takes too long.

\begin{figure}[htbp]
    \centering
    \includegraphics[width=12cm]{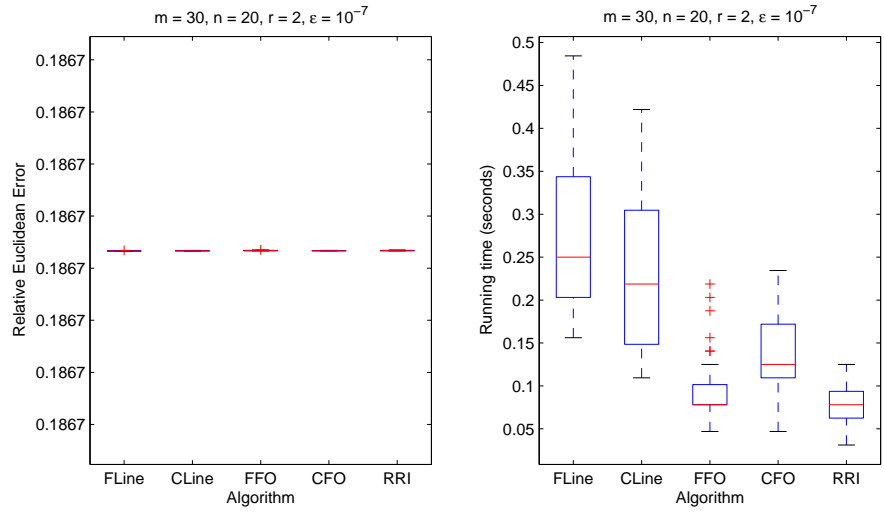}
    \caption{Comparison of selected algorithms for $\epsilon = 10^{-7}$}
    \label{fig:m30n20r2e7}
\end{figure}

\begin{figure}[htbp]
    \centering
    \includegraphics[width=12cm]{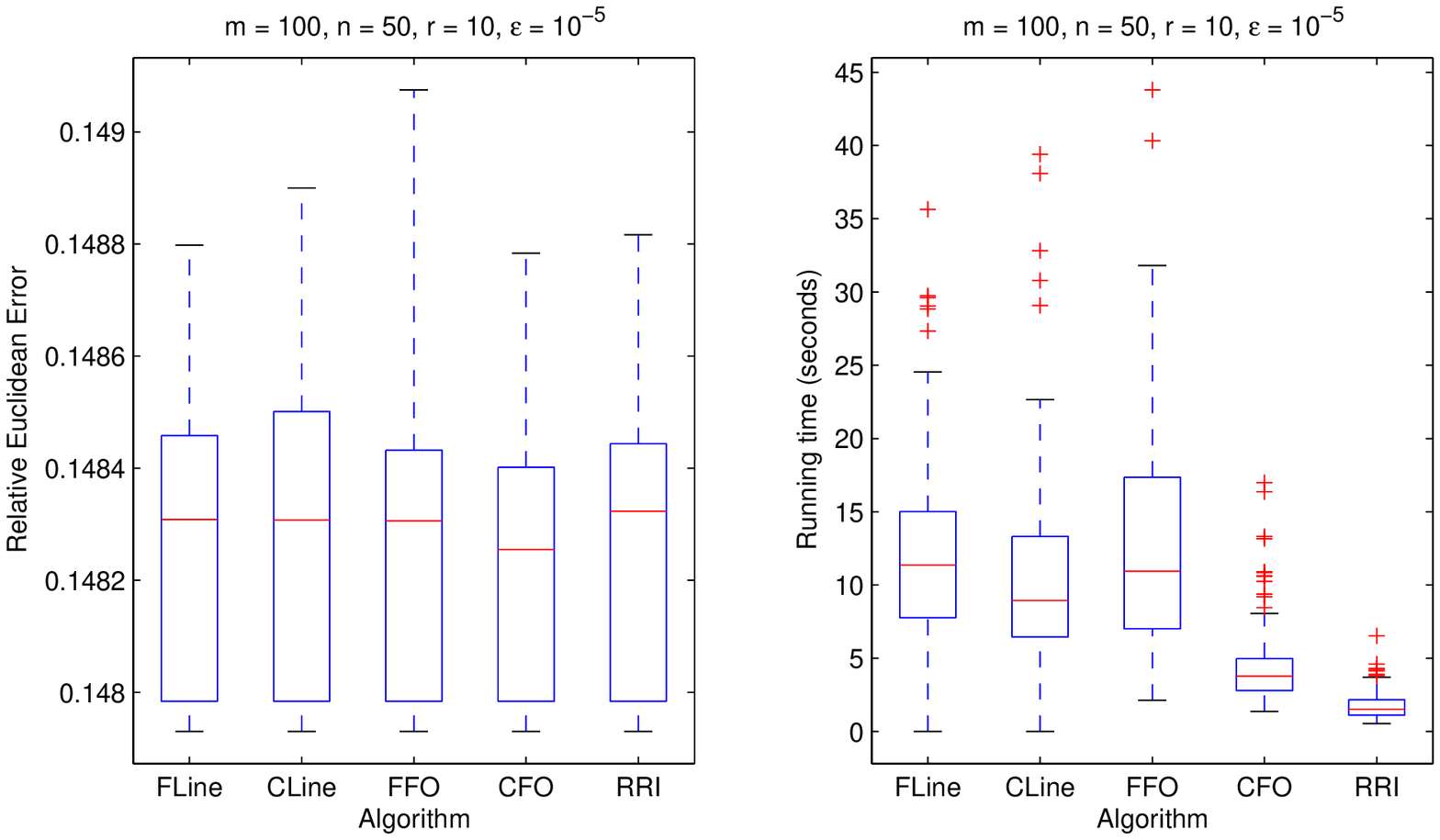}
    \caption{Comparison of selected algorithms for $\epsilon = 10^{-5}$}
    \label{fig:m100n50r10e5}
\end{figure}

\begin{figure}[htbp]
    \centering
    \includegraphics[width=12cm]{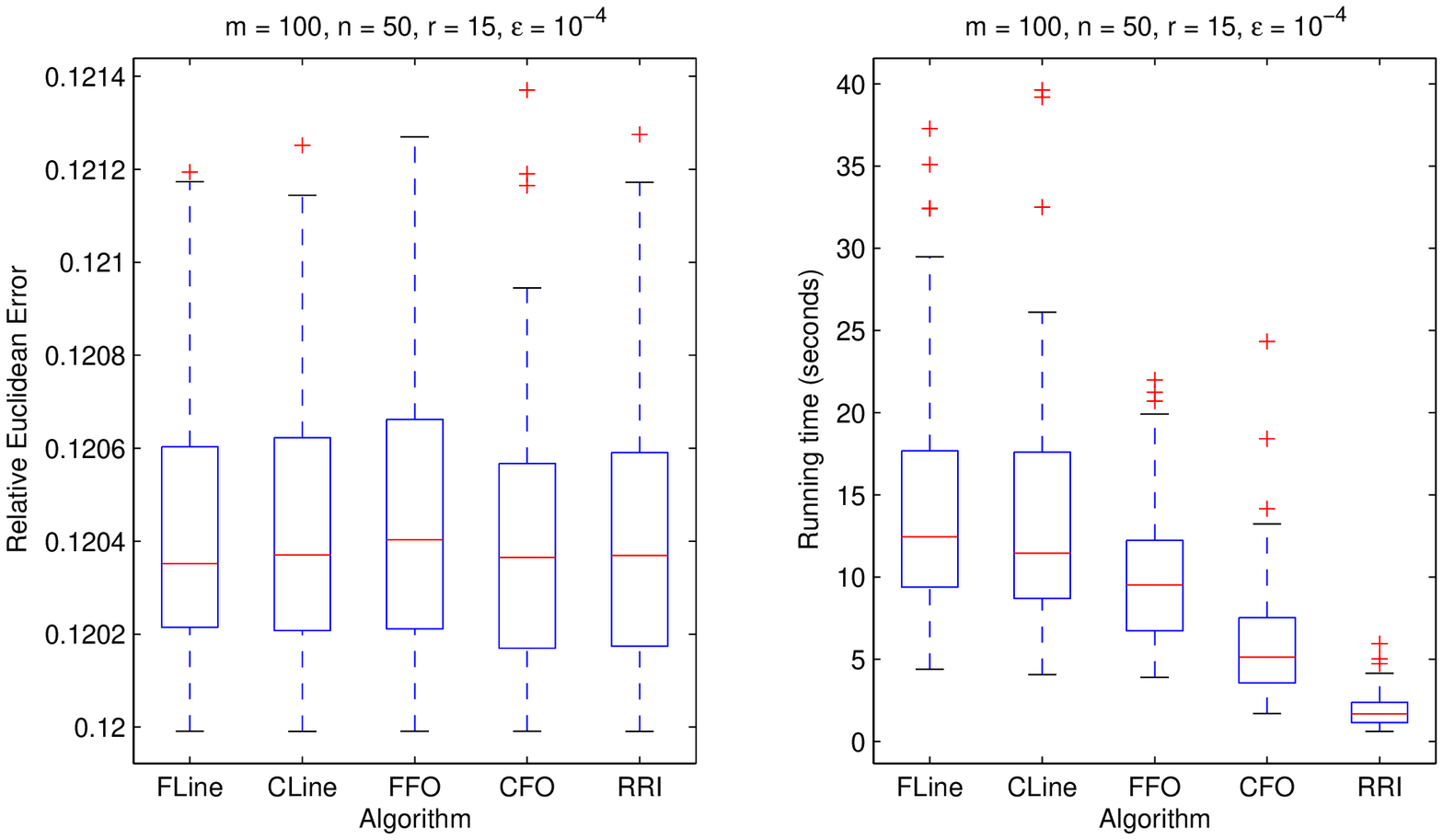}
    \caption{Comparison of selected algorithms for $\epsilon = 10^{-4}$}
    \label{fig:m100n50r15e4}
\end{figure}

\begin{figure}[htbp]
    \centering
    \includegraphics[width=12cm]{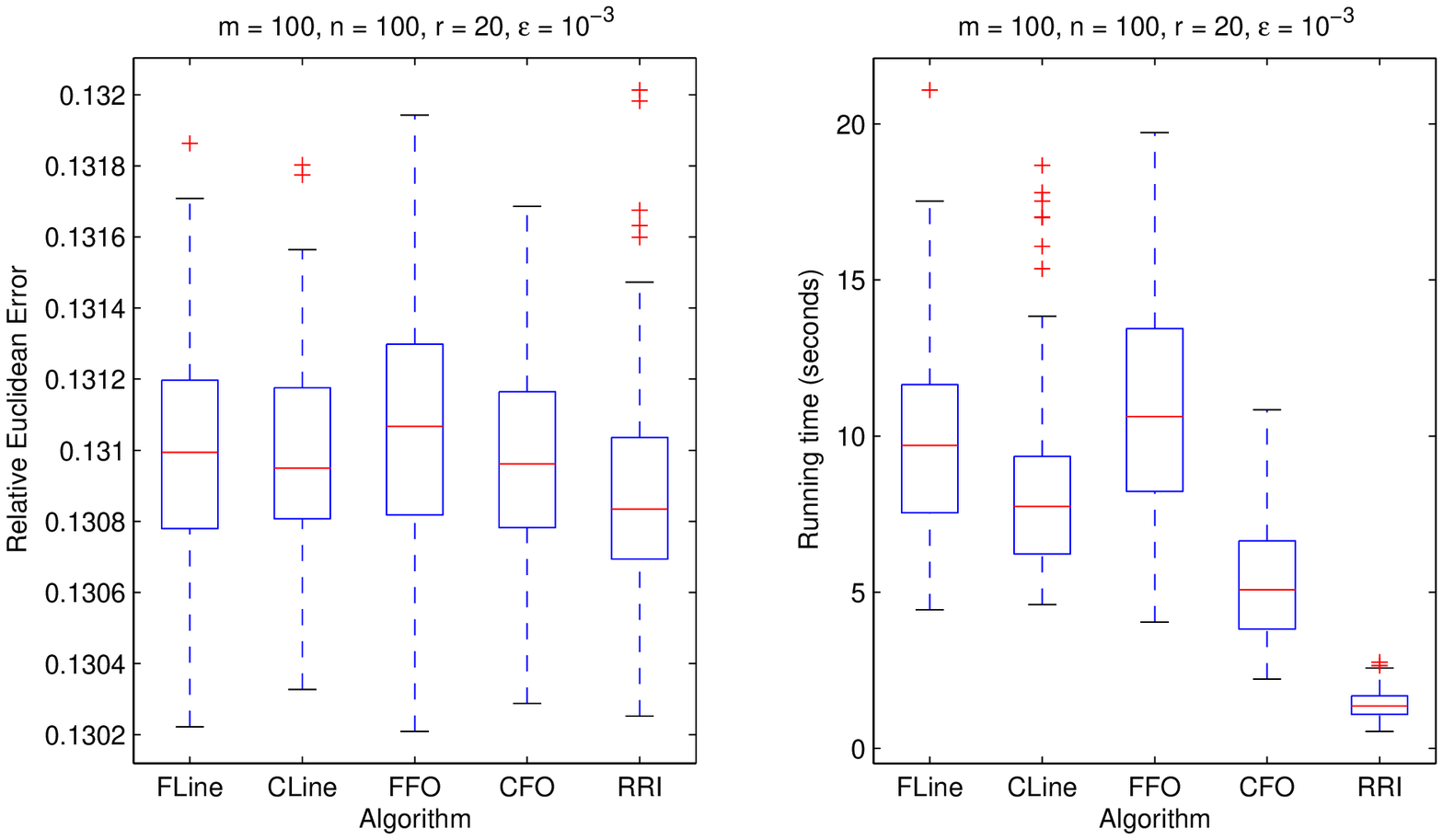}
    \caption{Comparison of selected algorithms for $\epsilon = 10^{-3}$}
    \label{fig:m100n100r20e3}
\end{figure}

We select five methods: \textbf{FLine}, \textbf{CLine}, \textbf{FFO}, \textbf{CFO} and \textbf{RRI} for a more detailed comparison. For each matrix $A$, we run these algorithms with $100$ different starting points. Figure \ref{fig:m30n20r2e7}, \ref{fig:m100n50r10e5}, \ref{fig:m100n50r15e4} and \ref{fig:m100n100r20e3} show the results with some different settings. One can see that, when the approximated errors are almost the same between the algorithms, \textbf{RRI} is the best overall in terms of running times.  It is probably because the RRI algorithm chooses only one vector $u_t$ or $v_t$ to optimize at once. This allows the algorithm to move \emph{optimally} down on partial direction rather than just a \emph{small step} on a more global direction. Furthermore, the computational load for an update is very small, only one matrix-vector multiplication is needed. All these factors make the running time of the RRI algorithm very attractive.

\subsection{Image data}

The following experiments use the Cambridge ORL face database as the input data. The database contains $400$ images of $40$ persons ($10$ images per person). The size of each image is $112 \times 92$ with $256$ gray levels per pixel representing a front view of the face of a person. The images are then transformed into 400 ``face vectors'' in $\mathbb R^{10304}$ ($112\times 92 = 10304$) to form the data matrix $A$ of size {$10304 \times 400$}. We used three weight matrices of the same size of $A$ (ie. {$10304 \times 400$}). Since it was used in \cite{leeseung99nature}, this data has become the standard benchmark for NMF algorithms.

\begin{figure}[htb]
    \centering
    \includegraphics[width=9cm]{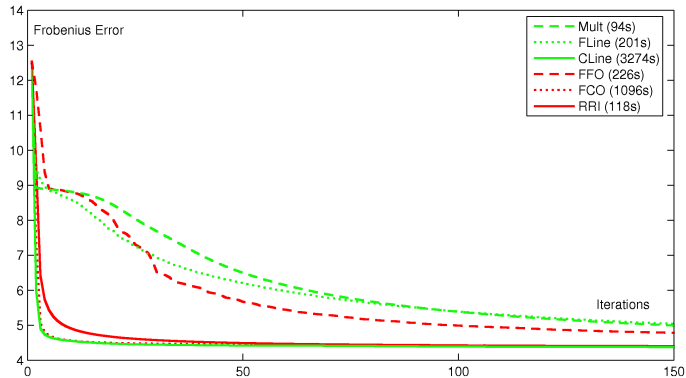}
    \caption{NMF: Error vs. Iterations}
    \label{fig:ErrvsIters}
\end{figure}

\begin{figure}[htbp]
    \centering
    \includegraphics[width=9cm]{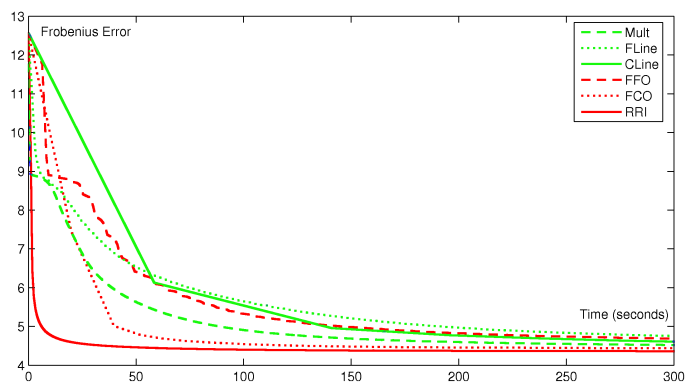}
    \caption{NMF: Error vs. Time}
    \label{fig:ErrvsTime}
\end{figure}

In the first experiment, we run six NMF algorithms described above on this data for the reduced rank of $49$. The original matrix $A$ is constituted by transforming each image into one of its column. Figure \ref{fig:ErrvsIters} shows for the six algorithms the evolution of the error versus the number of iterations. Because the minimization process is different in each algorithm, we will say that one iteration corresponds to all elements of both $U$ and $V$ being updated. Figure \ref{fig:ErrvsTime} shows the evolution of the error versus time. Since the work of one iteration varies from one algorithm to another, it is crucial to plot the error versus time to get a fair comparison between the different algorithms. In the two figures, we can see that the RRI algorithm behaves very well on this dataset. And since its computation load of each iteration is small and constant (without inner loop), this algorithm converges faster than the others.

\begin{figure}[htbp] 
    \centering
    \includegraphics[width=11cm]{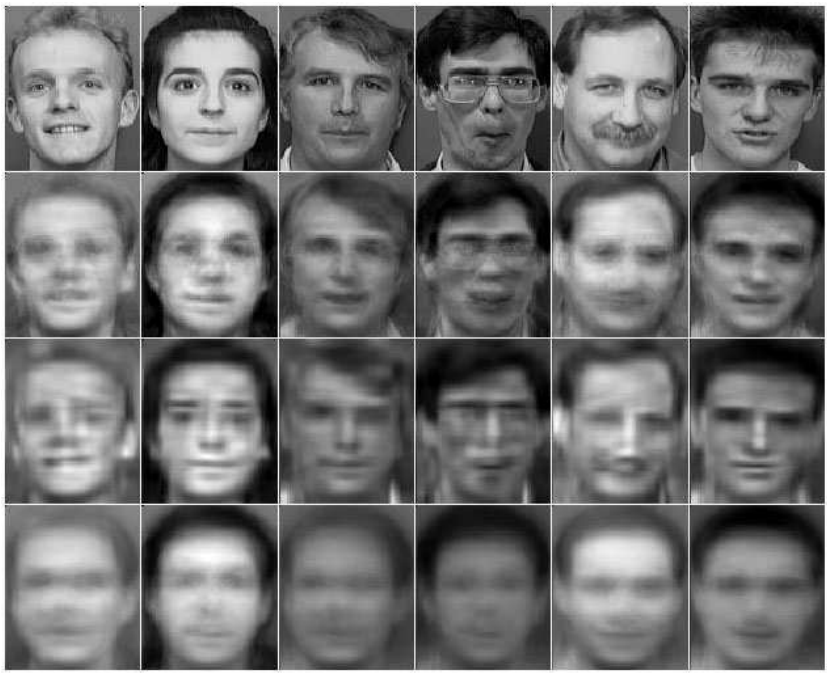}
    \caption{Tensor Factorization vs. Matrix Factorization on facial data. Six randomly chosen images from 400 of ORL dataset.  From top to bottom: original images, their $rank-8$ truncated SVD approximation, their $rank-142$ nonnegative tensor approximation ($150$ RRI iterations) and their $rank-8$ nonnegative matrix approximation ($150$ RRI iterations).}
    \label{fig:NTFFaces}
\end{figure}

In the second experiment, we construct a third-order nonnegative tensor approximation. We first build a tensor by stacking all $400$ images to have a $112 \times 92 \times 400$ nonnegative tensor. Using the proposed algorithm, a $rank-142$ nonnegative tensor is calculated to approximate this tensor. Figure \ref{fig:NTFFaces} shows the result for six images chosen randomly from the $400$ images. Their approximations given by the rank-142 nonnegative tensor are much better than that given by the rank-8 nonnegative matrix, even though they require similar storage space: $8*(112*92 + 400) = 85632$ and $142*(112 + 92 + 400) = 85768$. The rank-8 truncated SVD approximation (i.e. $[A_8]_+$) is also included for reference.

In the third experiment, we apply the variants of RRI algorithm mentioned in Section \ref{RRIVars} to the face databases. The following settings are compared:

\begin{itemize}
\item \textbf{Original}: original faces from the databases. 
\item \textbf{49NMF}: standard factorization (nonnegative vectors), $r = 49$. 
\item \textbf{100Binary}: columns of $U$ are limited to the scaled binary vectors, $r = 100$.
\item \textbf{49Sparse10}: columns of $U$ are sparse. Not more than $10\%$ of the elements of each column of $A$ are positive. $r = 49$.
\item \textbf{49Sparse20}: columns of $U$ are sparse. Not more than $20\%$ of the elements of each column of $A$ are positive. $r = 49$. 
\item \textbf{49HSparse60}: columns of $U$ are sparse. The Hoyer sparsity of each column of $U$ are $0.6$. $r = 49$.
\item \textbf{49HSparse70}: columns of $U$ are sparse. The Hoyer sparsity of each column of $U$ are $0.7$. $r = 49$. 
\item \textbf{49HBSparse60}: columns of $U$ are sparse. The Hoyer sparsity of each column of $U$ are $0.6$. Columns of $V$ are scaled binary. $r = 49$.
\item \textbf{49HBSparse70}: columns of $U$ are sparse. The Hoyer sparsity of each column of $U$ are $0.7$. Columns of $V$ are scaled binary. $r = 49$. 
\end{itemize}

\begin{figure}[htbp] 
    \centering
    \includegraphics[width=11cm]{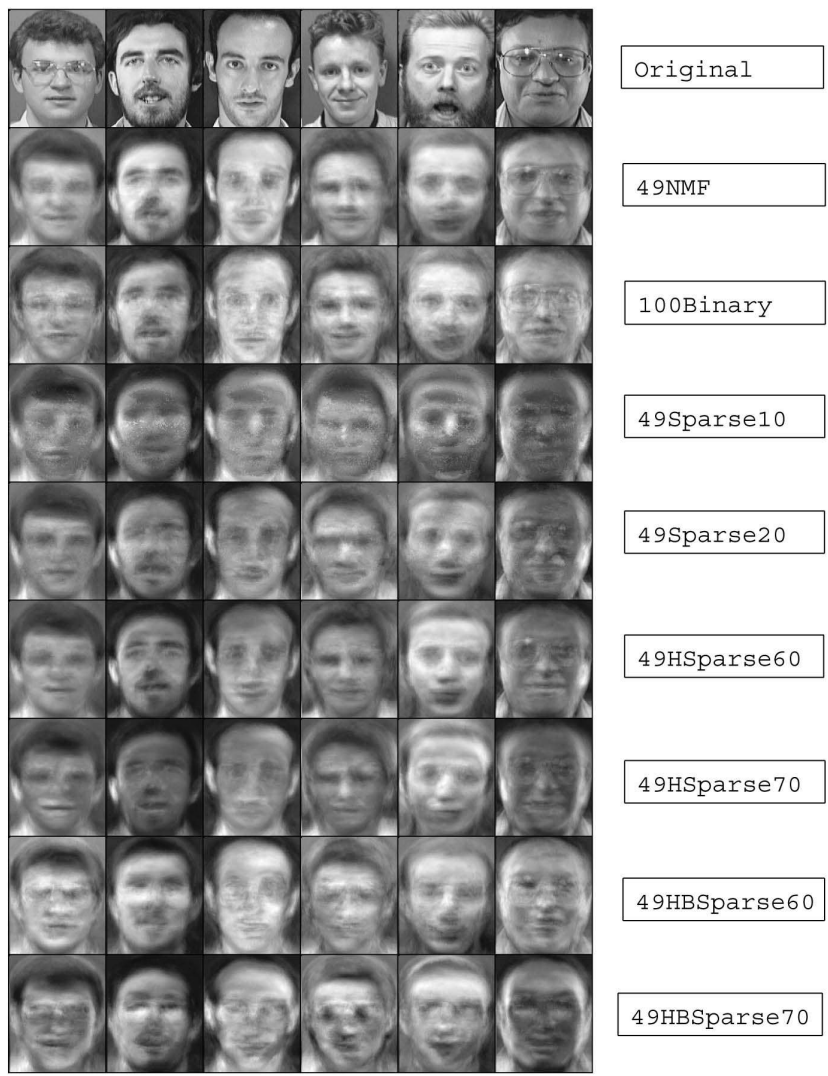}
    \caption{Nonnegative matrix factorization with several sparse settings}
    \label{fig:SparseFaces}
\end{figure}

\begin{figure}[htbp] 
    \centering
    \includegraphics[width=11cm]{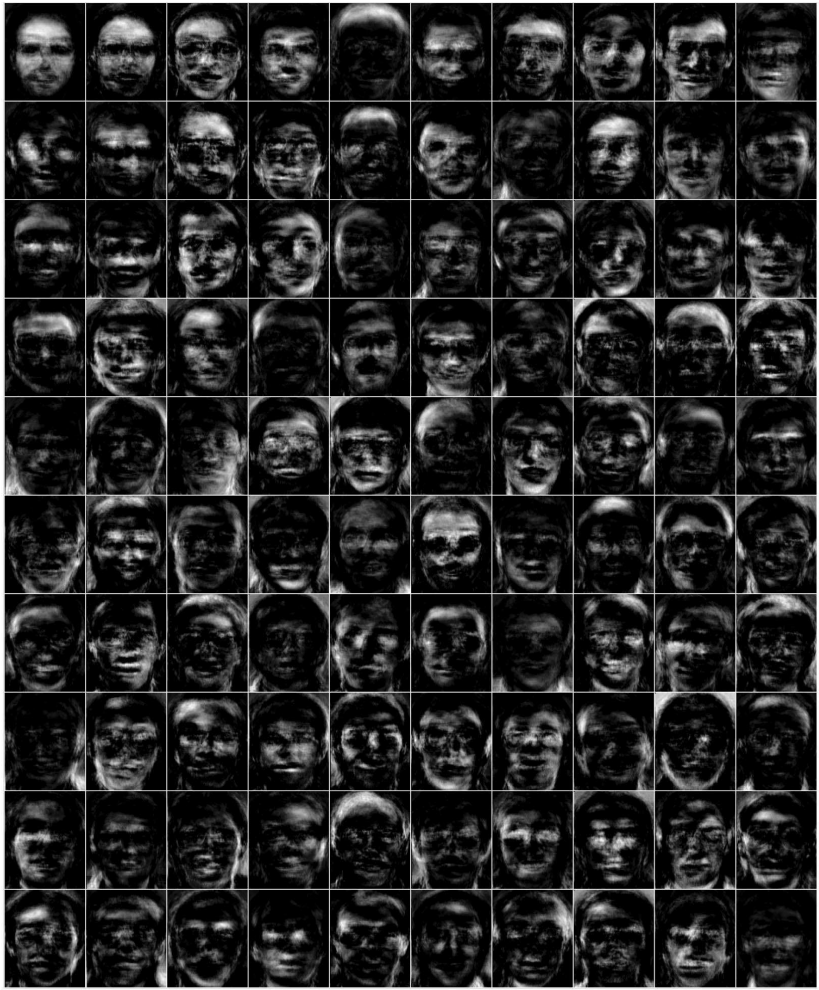}
    \caption{Bases from \textbf{100Binary} setting}
    \label{fig:Binary100}
\end{figure}

For each setting, we use RRI algorithm to compute the corresponding factorization. Some randomly selected faces are reconstructed by these settings as shown in Figure \ref{fig:SparseFaces}. For each setting, RRI algorithm produces a different set of bases to approximate the original faces. When the columns of $V$ are constrained to scaled binary vectors (\textbf{100Binary}), the factorization can be rewritten as $UV^T = \hat UB^T$, where $B$ is a binary matrix. This implies that each image is reconstructed by just the presence or absence of $100$ bases shown in Figure \ref{fig:Binary100}. 

\begin{figure}[htbp]
		\centering
		\begin{tabular}{cc}
				(a) & (b) \\
        \includegraphics[width=5.5cm]{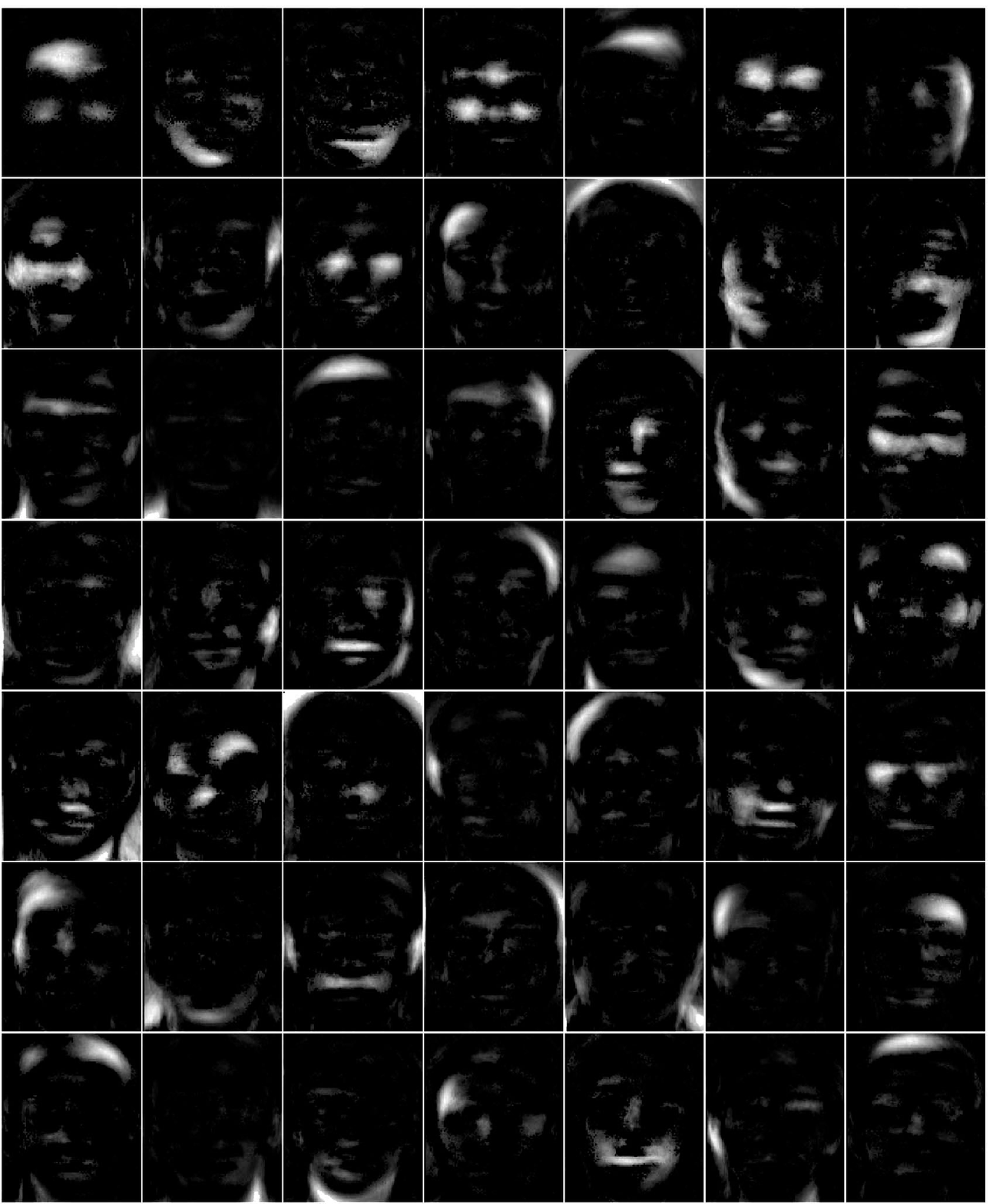}        
        &
        \includegraphics[width=5.5cm]{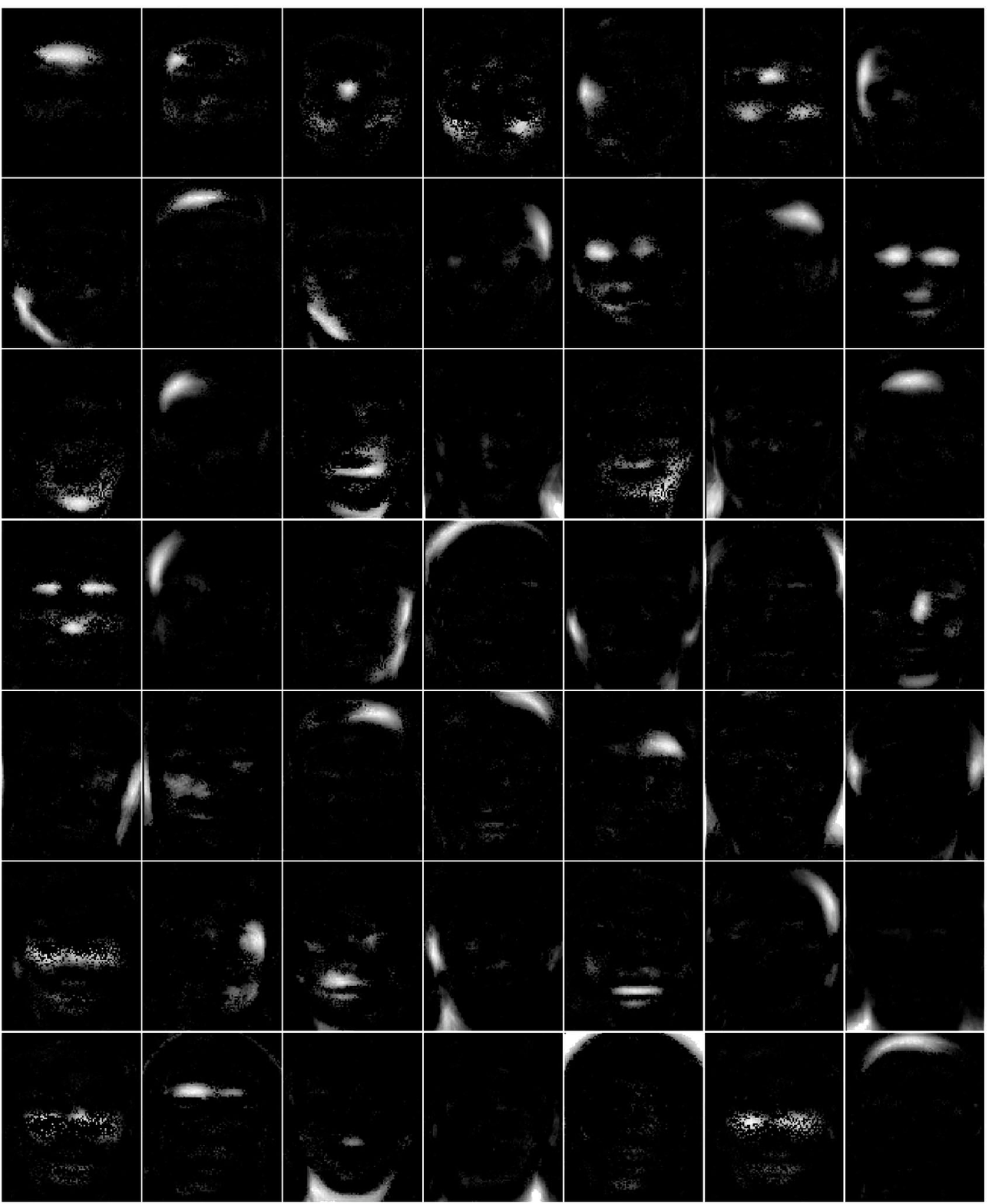}        
		\end{tabular}
		\caption{Sparse bases \textbf{49Sparse20} and \textbf{49Sparse10}. Maximal percentage of positive elements is $20\%$ (a) and $10\%$ (b)}.
		\label{fig:SparseBases}
\end{figure}

\begin{figure}[htbp]
		\centering
		\begin{tabular}{cc}
				(a) & (b) \\
        \includegraphics[width=5.5cm]{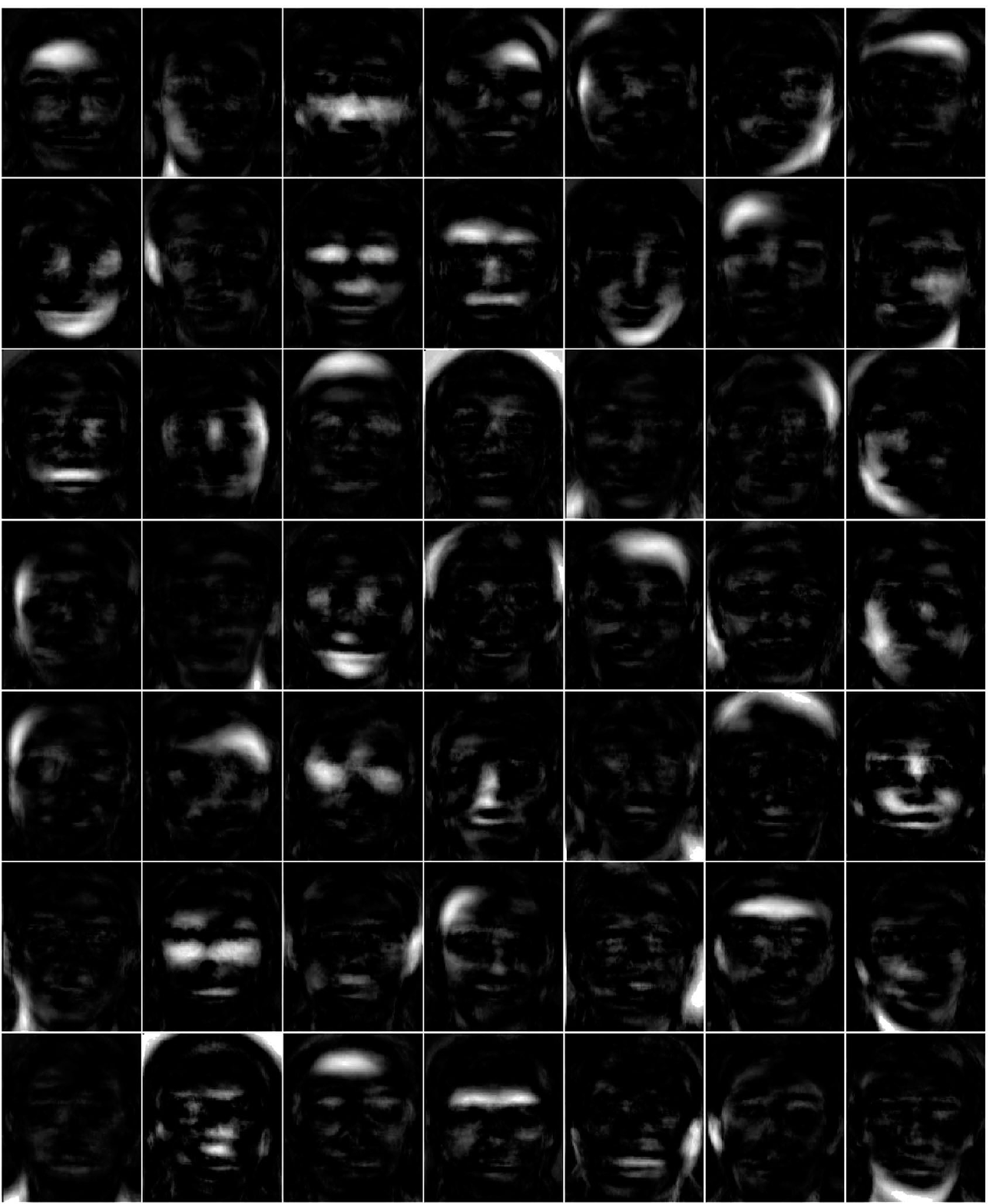}        
        &
        \includegraphics[width=5.5cm]{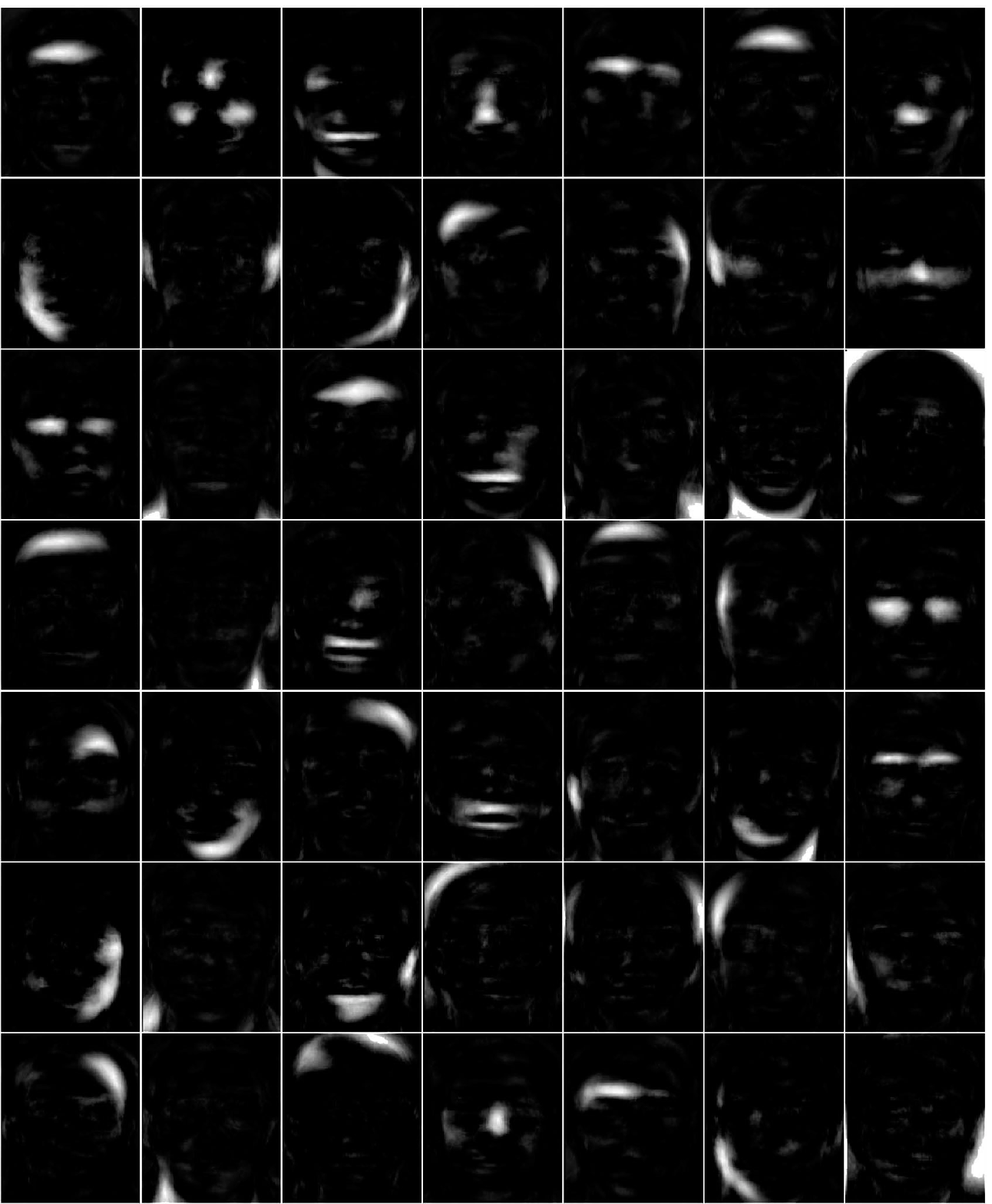}        
		\end{tabular}
		\caption{Hoyer sparse bases \textbf{49HSparse60} and \textbf{49HSparse70}. Sparsity of bases is $0.6$ (a) and $0.7$ (b)}.
		\label{fig:HSparseBases}
\end{figure}

Figure \ref{fig:SparseBases} and \ref{fig:HSparseBases} show nonnegative bases obtained by imposing some sparsity on the columns of $V$. The sparsity can be easily controlled by the percentages of positive elements or by the Hoyer sparsity measure. 

\begin{figure}[htbp]
		\centering
		\begin{tabular}{cc}
				(a) & (b) \\
        \includegraphics[width=5.5cm]{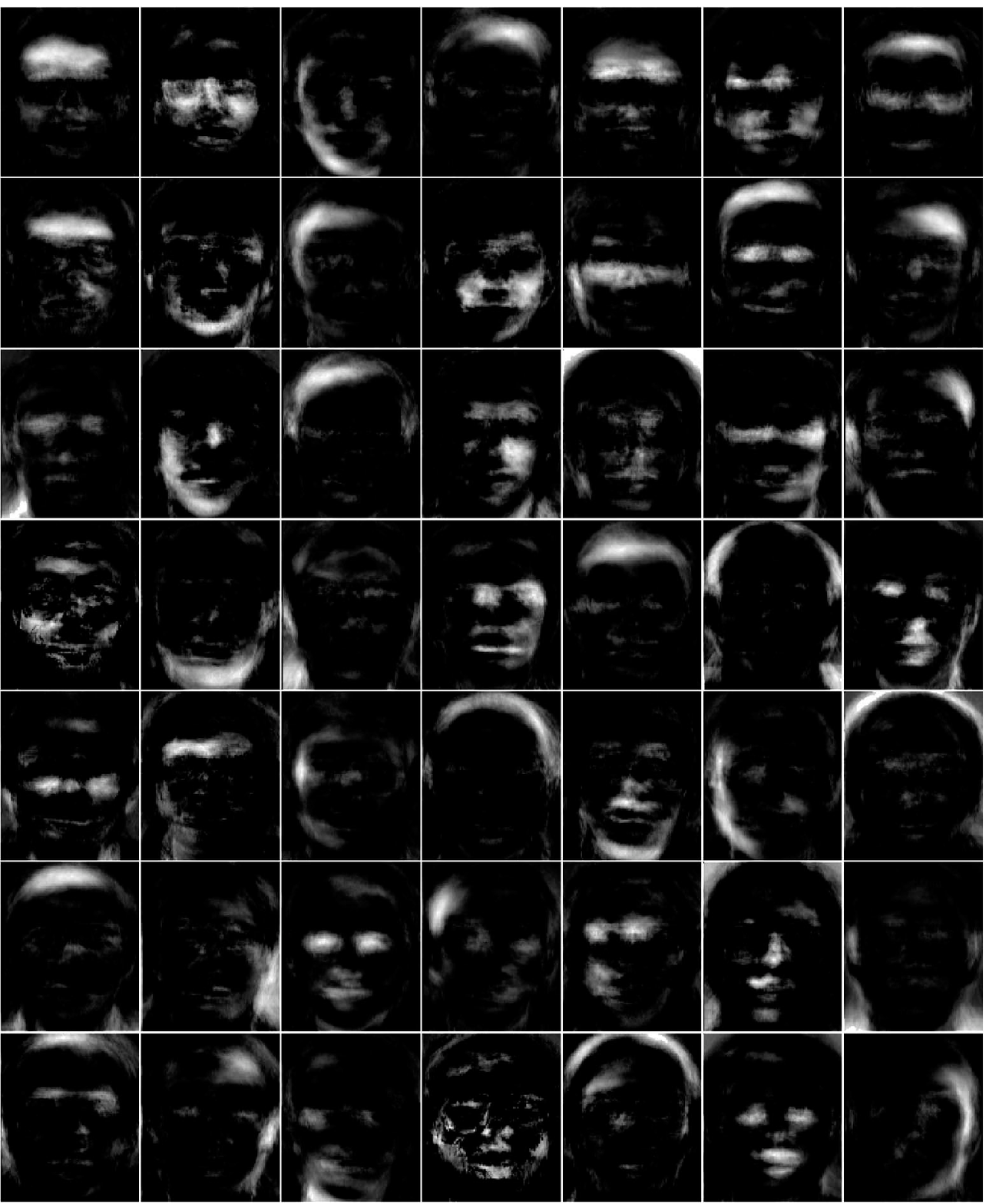}        
        &
        \includegraphics[width=5.5cm]{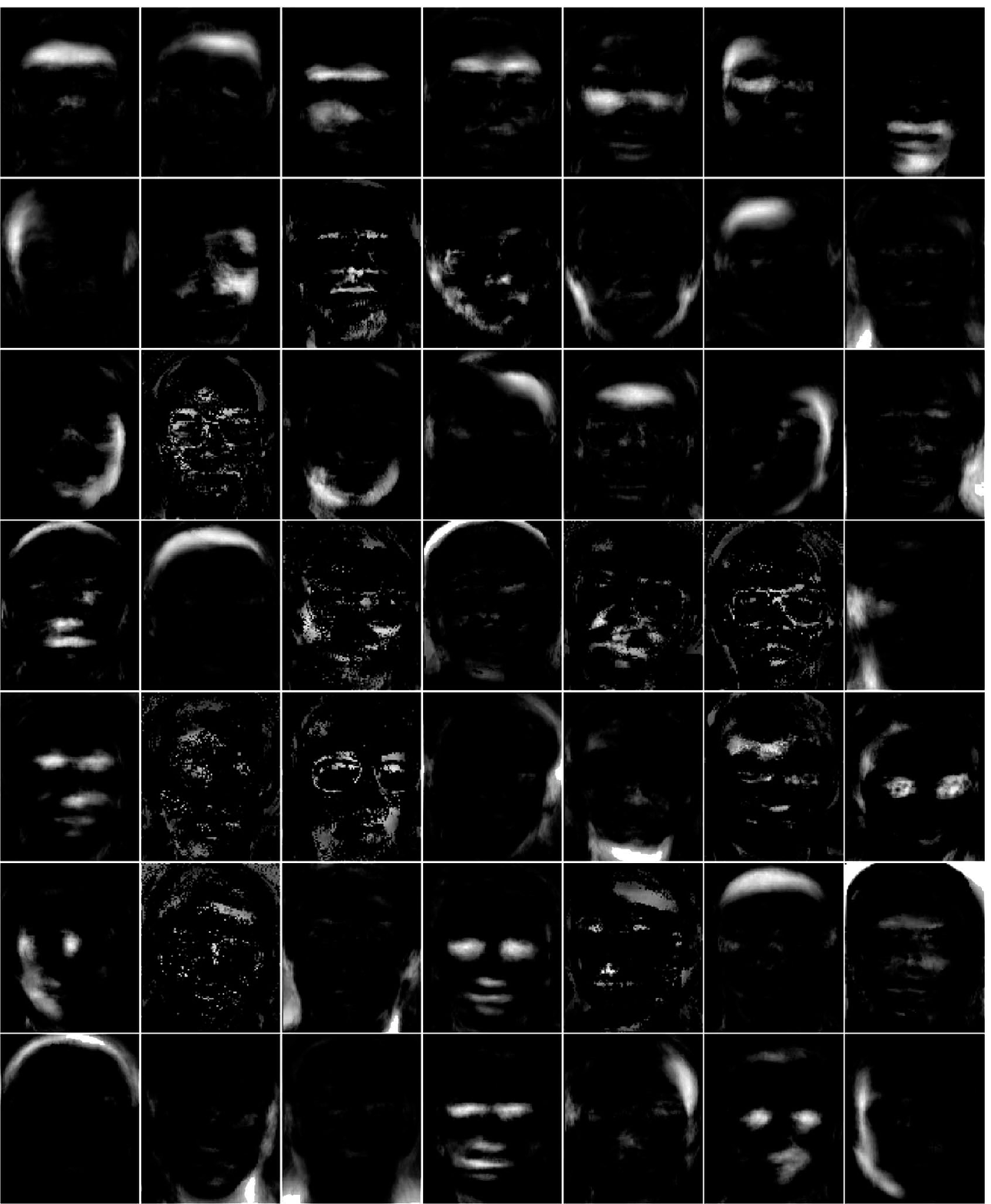}        
		\end{tabular}
		\caption{Hoyer sparse bases \textbf{49HBSparse60} and \textbf{49HBSparse70}. Sparsity of bases is $0.6$ (a) and $0.7$ (b). $V$ is binary matrix.}.
		\label{fig:HBSparseBases}
\end{figure}

Figure \ref{fig:HBSparseBases} combines the sparsity of the bases (columns of  $U$) and the binary representation of $V$. The sparsity is measured by the Hoyer measure as in Figure \ref{fig:HSparseBases}. Only with the absence or presence of these $49$ \emph{features}, faces are approximated as showed in the last two rows of Figure \ref{fig:SparseFaces}.

The above examples show how to use the variants of the RRI algorithm to control the sparsity of the bases. One can see that the sparser the bases are, the less storage is needed to store the approximation. Moreover, this provides a part-based decomposition using \emph{local} features of the faces. 

\subsection{Smooth approximation}

\begin{figure}[htbp] 
    \centering
    \includegraphics[width=11cm]{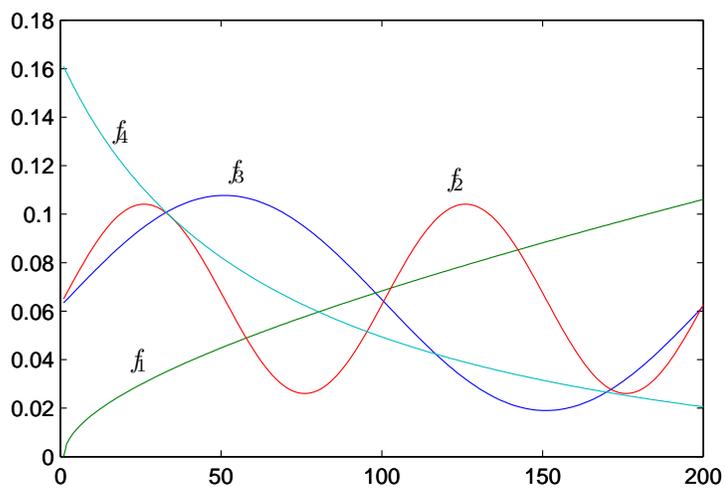}
    \caption{Smooth functions}
    \label{fig:smoothfunction}
\end{figure}

We carry out this experiment to test the new smoothness constraint introduced in the previous section:
$$
\frac{1}{2}\|R_i - u_iv^T\|^2_F + \frac{\delta}{2} \|v - B\hat v_i\|_F^2, \quad \delta > 0
$$
where $B$ is defined in (\ref{smoothmatrix}).

\begin{figure}[htb] 
    \centering
    \includegraphics[width=11cm]{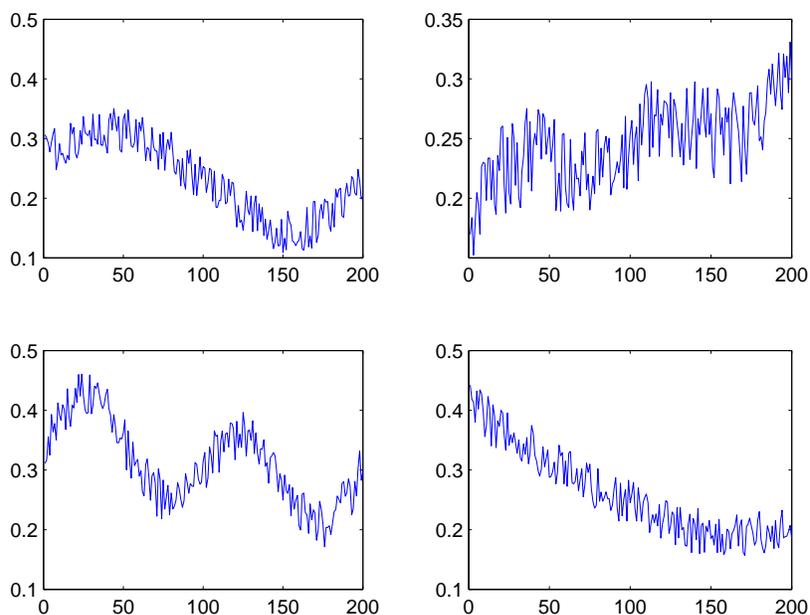}
    \caption{Randomly selected generated data}
    \label{fig:smoothmixture}
\end{figure}

We generate the data using four smooth nonnegative functions $f_1$, $f_2$, $f_3$ et $f_4$, described in Figure \ref{fig:smoothfunction}, where each function is represented as a nonnegative vector of size $200$.

\begin{figure}[htb] 
    \centering
    \includegraphics[width=11cm]{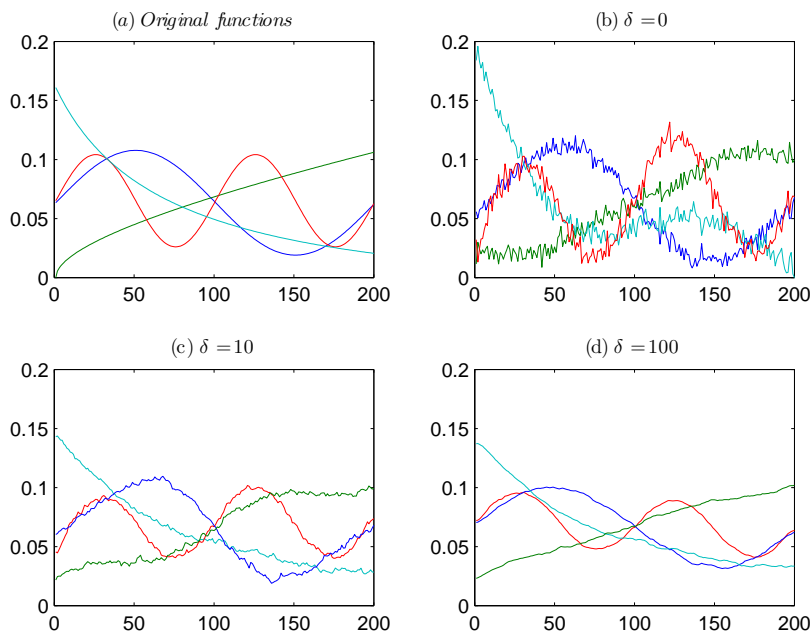}
    \caption{Original functions vs. reconstructed functions}
    \label{fig:smoothresult}
\end{figure}

We then generate a matrix $A$ containing $100$ mixture of these functions as follows
$$
A = max(FE^T + N, 0)
$$
where $F = [f_1 \ f_2 \ f_3 \ f_4]$, $E$ is a random nonnegative matrix and $N$ is normally distributed random noise with $\|N\|_F = 0.2 \|FE^T\|_F$. Four randomly selected columns of $A$ are plotted in Figure \ref{fig:smoothmixture}.

We run the regularized RRI algorithm to force the smoothness of  columns of $U$. We apply, for each run, the same value of $\delta$ for all the columns of $U$: $\delta = 0, 10, 100$. The results obtained through these runs are presented in Figure \ref{fig:smoothresult}. We see that, without regularization, i.e. $\delta = 0$, the noise is present in the approximation, which produces nonsmooth solutions. When increasing the regularizing terms, i.e. $\delta = 10, 100$, the reconstructed functions become smoother and the shape of the original functions are well preserved.  

This smoothing technique can be used for applications like that in \cite{pauca2006nmf}, where smooth spectral reflectance data from space objects is unmixed. The multiplicative rules are modified by adding the two-norm regularizations on the factor $U$ and $V$ to enforce the smoothness. This is a different approach, therefore, a comparison should be carried out.

We have described a new method for nonnegative matrix factorization that has a good and fast convergence. Moreover, it is also very flexible to create variants and to add some constraints as well. The numerical experiments show that this method and its derived variants behave very well with different types of data. This gives enough motivations to extend to other types of data and applications in the future. In the last two chapters of this thesis, it is applied to weighted cost functions and to symmetric factorizations. 

\section{Conclusion}

This paper focuses on the descent methods for Nonnegative Matrix Factorization, which are characterized by nonincreasing updates at each iteration. 

We present also the Rank-one Residue Iteration algorithm for computing an approximate Nonnegative Matrix Factorization. It uses recursively nonnegative rank one approximations of a residual matrix that is not necessarily nonnegative. This algorithm requires no parameter tuning, has nice properties and typically converges quite fast. It also has many potential extensions. During the revision of this report, we were informed that essentially the same algorithm was published in an independent contribution \cite{Cichocki2007HALS} and also mentioned later in an independent personal communication \cite{gillis2007}. 

\section*{Acknowledgments}
This paper presents research results of the Concerted Research Action(ARC) "Large Graphs and Networks" of the French Community of Belgium and the Belgian Network DYSCO (Dynamical Systems, Control, and Optimization), funded by the Interuniversity Attraction Poles Programme, initiated by the Belgian State,  Science Policy Office. The scientific responsibility rests with its authors. Ngoc-Diep Ho is a FRIA fellow.



%




\end{document}